
\input amstex
\documentstyle{amsppt}

\hcorrection{19mm}

\nologo
\NoBlackBoxes


\topmatter 

\title   The classification of Toroidal Dehn surgeries on Montesinos knots
\endtitle
\author  Ying-Qing Wu$^1$
\endauthor
\leftheadtext{Ying-Qing Wu}
\rightheadtext{Toroidal surgery on arborescent knots}
\address Department of Mathematics, University of Iowa, Iowa City, IA
52242
\endaddress
\email  wu\@math.uiowa.edu
\endemail

\thanks  $^1$ Partially supported by NSF grant \#DMS 0203394
\endthanks

\abstract Exceptional Dehn surgeries have been classified for 2-bridge
knots and Montesinos knots of length at least 4.  In this paper we
classify all toroidal Dehn surgeries on Montesinos knots of length 3.
\endabstract

\endtopmatter

\document

\define\proof{\demo{Proof}}
\define\endproof{\qed \enddemo}
\define\a{\alpha}
\define\bu{{\bar u}}
\redefine\b{\beta}

\redefine\r{\gamma}

\redefine\bdd{\partial}
\redefine\cd{{\Cal D}}
\define\Int{\text{\rm Int}}

\define\lcm{\text{\rm lcm}}

\define\<{\langle}
\define\>{\rangle}

\TagsOnRight

\input epsf.tex

\head 1. Introduction \endhead

A Dehn surgery on a hyperbolic knot $K$ along a slope $\delta$ is said
to be {\it exceptional\/} if the resulting manifold $K_{\delta}$ is
either reducible, toroidal, or a small Seifert fibered manifold.  If
the Geometrization Conjecture [Th] is true then these are exactly the
surgeries such that $K_{\delta}$ is non-hyperbolic.  By Thurston's
Hyperbolic Surgery Theorem, all but finitely many Dehn surgeries on a
hyperbolic knot produce hyperbolic manifolds, hence there are only
finitely many exceptional surgeries.

It is known that there are no exceptional surgeries on Montesinos
knots of length at least four [Wu1].  Exceptional surgeries for
2-bridge knots have been classified in [BW].  Thus length 3 knots are
the only ones among the Montesinos knots which have not been settled.
In this paper we will classify toroidal surgeries for such knots.  See
Theorems 1.1 and 1.2 below.  By [Eu1] there is no reducible surgery on
a hyperbolic Montesinos knot because it is strongly invertible.  It
remains a challenging open problem to determine all small Seifert
fibred surgeries on Montesinos knots of length 3.

Hatcher and Oertel [HO] have an algorithm to determine all boundary
slopes of a given Montesinos knot.  We will therefore focus on finding
all length 3 knots such that some of their boundary slopes are
toroidal slopes.  Each incompressible surface in the exterior of $K$
corresponds to three ``allowable edgepaths'' $\r_1, \r_2, \r_3$.  We
will define an Euler number for allowable edgepaths, and show that if
$F(\r_1, \r_2, \r_3)$ is a punctured torus then one of the $\r_i$ must
have non-negative Euler number.  We then analyze the graph of
Hatcher-Oertel (Figure 2.1), and show that the ending point of one of
the above $\r_i$ must lie in a subgraph consisting of 7 edges.  This
breaks the problem down to several different cases.  We will then use
the properties of allowable edgepaths to find all possible solutions
in each case.

Two knots are considered {\it equivalent\/} if there is a (possibly
orientation reversing) homeomorphism of $S^3$ sending one knot to
the other.  Thus $K_1$ is equivalent to $K_2$ if $K_1$ is isotopic to
$K_2$ or its mirror image.  Similarly, if $N(K_i)$ is a neighborhood
of $K_i$ and $\delta_i$ is a slope on $\bdd N(K_i)$, then $(K_1,
\delta_1)$ is equivalent to $(K_2, \delta_2)$ if there is a
homeomorphism of $S^3$ sending $N(K_1)$ to $N(K_2)$ and $\delta_1$ to
$\delta_2$.  The following is the classification theorem for toroidal
boundary slopes of Montesinos knots of length 3.  Some knots are
listed more than once, with different boundary slopes, which means
that they admit more than one toroidal surgery.  The variable $\bu$
is the $u$ coordinate of the ending points of the edgepaths, which
will be defined in Section 2.  Note that some knots may have the same
toroidal boundary slope at different $\bu$ values, in which case we
will only list one $\bu$ value.

\proclaim{Theorem 1.1} Let $K$ be a hyperbolic Montesinos knot of
length $3$, let $E(K) = S^3 - \Int N(K)$, and let $\delta$ be a slope
on $\bdd E(K)$.  Then $E(K)$ contains an essential surface $F$ of
genus one with boundary slope $\delta$ if and only if $(K, \delta)$ is
equivalent to one of the pairs in the following list.

(1) $K = K(1/q_1,\; 1/q_2,\; 1/q_3)$, $q_i$ odd, $|q_i|>1$,
    $\delta = 0$; $\bu = 1$.

(2) $K = K(1/q_1,\; 1/q_2,\; 1/q_3)$, $q_1$ even, $q_2, q_3$ odd,
 $|q_i|>1$, $\delta = 2(q_2+q_3)$; $\bu = 1$.

(3) $K = K(-1/2,\; 1/3,\; 1/(6+1/n))$, $n \neq 0, -1$,  $\delta = 16$
if $n$ is odd, and $0$ if $n$ is even; $\bu = 6$.

(4) $K = K(-1/3,\; -1/(3+1/n),\; 2/3)$, $n\neq 0, -1$, $\delta = -12$
when $n$ is odd, and $\delta=4$ when $n$ is even; $\bu = 3$.

(5) $K = K(-1/2,\; 1/5,\; 1/(3+1/n))$, $n$ even, and $n\neq 0$,
$\delta = 5 - 2n$; $\bu = 3$.

(6) $K = K(-1/2,\; 1/3,\; 1/(5+1/n))$, $n$ even, and $n\neq 0$,
$\delta = 1-2n$;  $\bu = 3$.

(7) $K = K(-1/(2+1/n),\; 1/3,\; 1/3)$, $n$ odd, $n \neq -1$,
$\delta = 2n$;  $\bu = 2$.

(8) $K = K(-1/2,\; 1/3,\; 1/(3+1/n))$, $n$ even, $n \neq 0$,
$\delta = 2-2n$;  $\bu = 2$.

(9) $K = K(-1/2,\; 2/5,\; 1/9)$, $\delta = 15$;  $\bu = 5$.

(10) $K = K(-1/2,\; 2/5,\; 1/7)$, $\delta = 12$;  $\bu = 4$.

(11) $K = K(-1/2,\; 1/3,\; 1/7)$, $\delta = 37/2$;  $\bu = 2.5$.

(12) $K = K(-2/3,\; 1/3,\; 1/4)$, $\delta = 13$;   $\bu = 2.5$.

(13) $K = K(-1/3,\; 1/3,\; 1/7)$, $\delta = 1$;  $\bu = 2.5$.
\endproclaim

For each case in Theorem 1.1, the candidate system $(\r_1, \r_2,
\r_3)$ is given in the proofs of the lemmas, hence it is straight
forward using the algorithm of Hatcher-Oertel to calculate the
boundary slope of $F(\r_1, \r_2, \r_3)$ and show that it is an
incompressible toroidal surface.  For each individual knot this can
also be verified using a computer program of Dunfield [Dn].  We will
therefore concentrate on showing the ``only if'' part, that is, if the
exterior of $K$ has an incompressible toroidal surface with boundary
slope $\delta$ then $(K, \delta)$ must be one of those in the list.

In general, the existence of a toroidal incompressible surface $F$
with boundary slope $\delta$ in the exterior of a knot $K$ does not
guarantee that $K_{\delta}$ is toroidal, because the corresponding
closed surface $\hat F$ may be compressible in $K_{\delta}$.  However,
the following theorem shows that this does not happen for Montesinos
knots of length 3; hence the above theorem actually gives a
classification of all toroidal surgeries for Montesinos knots of
length 3.

\proclaim{Theorem 1.2} Let $K$ be a Montesinos knot of length 3, and
let $\delta$ be a slope on $T = \bdd N(K)$.  Then $K_{\delta}$ is
toroidal if and only if $(K, \delta)$ is equivalent to one of those in
the list of Theorem 1.1.
\endproclaim

Together with [BW] and [Wu1], this gives a complete classification of
toroidal surgeries on all Montesinos knots.  The following are some of
the consequences.

(1) The only non-integral toroidal surgery on a Montesinos knot is the
    $37/2$ surgery on $K(-1/2,\; 1/3,\; 1/7)$.

(2) No Montesinos knot admits more than three toroidal surgeries, and the
    Figure 8 knot and $K(-1/2,\; 1/3,\; 1/7)$ are the only ones
    admitting three toroidal surgeries.

(3) By [BW], a 2-bridge knot admits exactly two toroidal surgeries if
    and only if it is associated to the rational number $1/(2+1/n)$
    for some $|n|>2$.  By checking the list in Theorem 1.1 for knots
    which are listed more than once, we see that $K(t_1, t_2, t_3)$
    admits exactly two toroidal surgeries if and only if it is
    equivalent to one of the following 5 knots.
$$
\align
& K(-1/2,\; 1/3,\; 2/11), \qquad \text{$\delta = 0$ and $-3$;} \\
& K(-1/3,\; 1/3,\; 1/3), \qquad \text{$\delta = 0$ and $2$;} \\
& K(-1/3,\; 1/3,\; 1/7), \qquad \text{$\delta = 0$ and $1$;} \\
& K(-2/3,\; 1/3,\; 1/4), \qquad \text{$\delta = 12$ and $13$;} \\
& K(-1/3,\; -2/5,\; 2/3), \qquad \text{$\delta = 4$ and $6$.} 
\endalign
$$

(4) A toroidal essential surface $F$ in Theorem 1.1 has at most 4
boundary components.  In case (1) of Theorem 1.1 $F$ is a Seifert
surface with a single boundary component.  In all other cases $F$ is a
separating surface and the result follows from the proof of Theorem
1.2.

\medskip

(1) also follows from results of Gordon and Luecke [GL2] and
Eudave-Munoz [Eu2], which classified non-integral toroidal surgeries
on all knots in $S^3$.  There are many other interesting results about
toroidal Dehn surgery, see for example [Go, GL1, GL3, GW, Oh, Te1,
Te2, Wu2].

The paper is organized as follows.  In Section 2 we give a brief
introduction to some definitions and results of Hatcher and Oertel in
[HO], then define and explore the properties of Euler numbers $e(\r)$
for any edgepath in the Hatcher-Oertel graph $\cd$ shown in Figure 2.1.
It will be shown that if $F(\r_1, \r_2, \r_3)$ is a punctured torus
then up to equivalence the ending point $v_1$ of $\r_1$ must lie on
the subgraph $G$ in Figure 2.4.  Sections 3, 4 and 5 discuss the
cases that $v_1$ lies on a horizontal edge in $G$, and Section 6 deals
with the remaining cases.  The proofs of Theorems 1.1 and 1.2 will be
given in Section 7.

I would like to thank Hatcher and Oertel for their algorithm in [HO],
which is crucial to the current work.  Thanks also to Nathan Dunfield
for his program [Dn], which has been used to verify that the slopes
found in Theorem 1.1 are indeed toroidal boundary slopes.

\head 2.  Preliminaries \endhead

In this section we first give a brief introduction to some results of
Hatcher-Oertel in [HO].  We will then define Euler numbers for points
and edgepaths on the Hatcher-Oertel diagram $\cd$, and show how they
are related to the Euler characteristic of the corresponding surfaces.
The main result is Theorem 2.8, which will play a key role in finding
Montesinos knots which admit toroidal surgeries.

\subhead  2.1.  The diagram $\cd$  \endsubhead

The diagram $\cd$ of Hatcher-Oertel is a 2-complex on the plane $\Bbb
R^2$ consisting of vertices, edges and triangular faces described as
follows.  See Figure 2.1, which is the same as [HO, Figure 1.3].
Unless otherwise stated, we will always write a rational number as
$p/q$, where $p,q$ are coprime integers, and $q>0$.

\roster
\item To each rational number $y = p/q$ is associated a vertex
    $\< y \>$ in $\cd$, which has Euclidean coordinates $(x,y) =
    ((q-1)/q,\; p/q)$.

\item For each rational number $y= p/q$, there is also an ``ideal''
    vertex $\< y \>_0$ with Cartesian coordinates $(1, p/q)$.

\item There is a vertex $\infty = 1/0$ located at $(-1, 0)$.

\item There is an edge $E = \< p/q,\; p'/q' \>$ connecting $\< p/q \>$ to
    $\< p'/q' \>$ if and only if $|pq'-p'q| = 1$.  Thus for example,
    there is an edge connecting $\infty$ to each vertex $\< p/1 \>$,
    and there is an edge connecting $\< p/q \>$ to $\< 0 \>$ if and
    only if $p=\pm 1$.

\item For each rational number $y$ there is a horizontal edge
    $L(y)$ connecting $\< y \>$ to the ideal vertex $\< y \>_0$.

\item A face of $\cd$ is a triangle bounded by three
  non-horizontal edges of $\cd$.
\endroster

\bigskip
\leavevmode

\centerline{\epsfbox{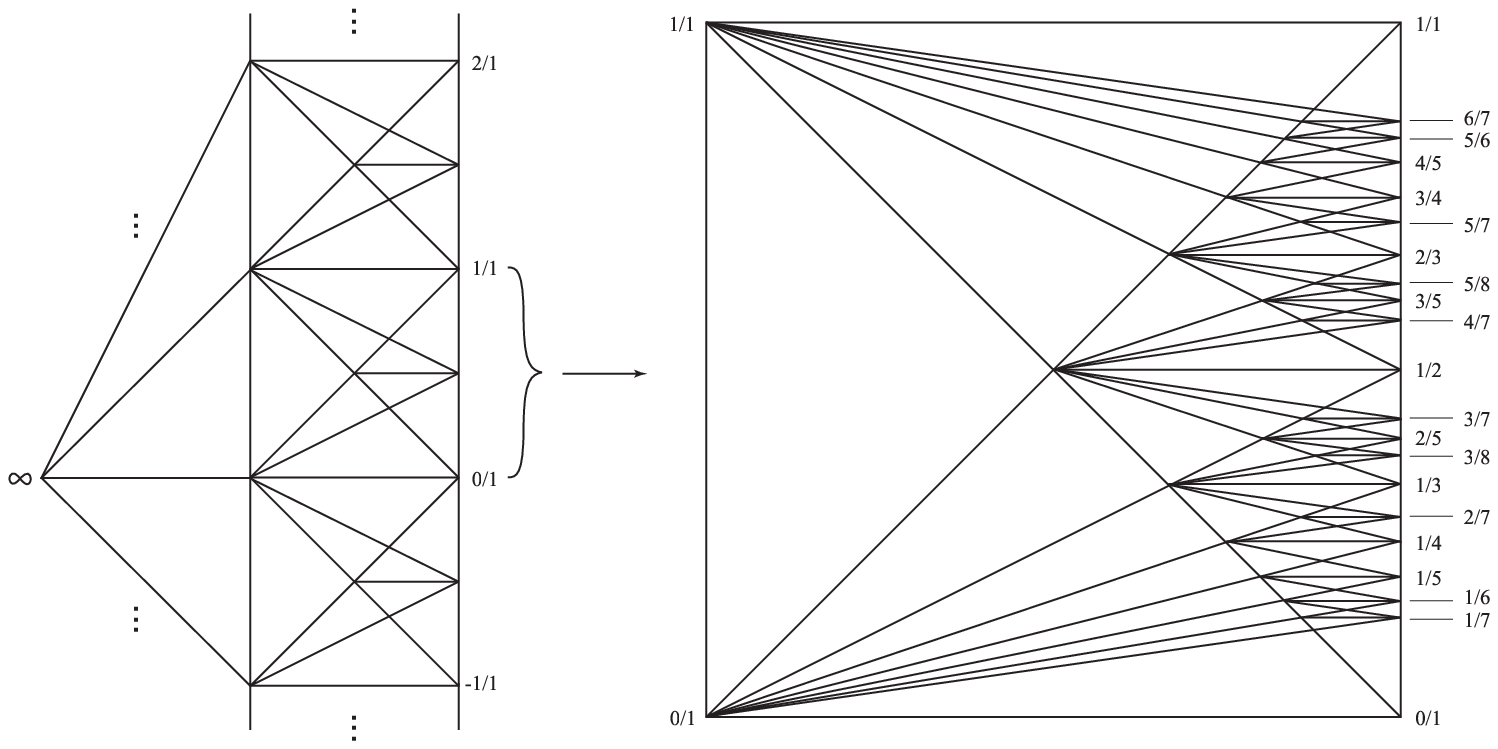}}
\bigskip
\centerline{Figure 2.1}
\bigskip

Note that a non-horizontal line segment in the figure with one
endpoint on an ideal vertex (i.e.\ a vertex on the vertical line
$x=1$) is not an edge of $\cd$.  It is a union of infinitely many
edges of $\cd$ and contains infinitely many vertices of $\cd$.
Similarly a triangle $\Delta$ in the figure is not a face of
$\cd$ if it contains a horizontal edge because there are edges in the
interior of $\Delta$.  Actually in this case $\Delta$ is a union of
infinitely many faces of $\cd$.  On the other hand, if all three
vertices of a triangle $\Delta$ in the figure are non-ideal vertices
of $\cd$, and if all three boundary edges of $\Delta$ are edges of
$\cd$ as defined above, then $\Delta$ is a face of $\cd$; in
particular, its interior contains no other edges or vertices of $\cd$.

\subhead 2.2.  Allowable edgepaths, candidate systems, and candidate
surfaces 
\endsubhead

An {\it edgepath\/} $\gamma$ in $\cd$ is a piecewise linear path in
the 1-skeleton of $\cd$.  Note that the endpoints of $\gamma$ may not
be vertices of $\cd$.  An edgepath $\r$ is a {\it constant path\/} if
its image is a single point.

Let $K = K(t_1,\; t_2,\; t_3)$ be a Montesinos knot of length 3.  Let
$\r_1, \r_2, \r_3$ be three edgepaths in $\cd$.  According to [HO,
P457], we say that the three edgepaths form a {\it candidate system\/}
for $K(t_1, t_2, t_3)$ if they satisfy the following conditions.

\roster
\item The starting point of $\r_i$ is on the horizontal edge $L(t_i)$, and
  if this starting point is not the vertex $\< t_i \>$ then $\r_i$ is
  a constant path.

\item  $\r_i$ is minimal in the sense that it never stops and
  retraces itself, or goes along two sides of a triangle of $\cd$ in
  succession.

\item  The ending points of $\r_i$ are rational points
  $\cd$ which all lie on one vertical line and whose vertical
  coordinates add up to zero.

\item  $\r_i$ proceeds monotonically from right to left,
  ``monotonically'' in the weak sense that motion along vertical edges
  is permitted.
\endroster

Each $\r_i$ above is called an {\it allowable edgepath.}  By
definition, an allowable edgepath must be of one of the following
three types.

\roster
\item  A constant path on a horizontal edge, possibly at a vertex of
  $\cd$.  

\item  An edgepath with both endpoints on vertices of $\cd$.

\item  An edgepath starting from a vertex of $\cd$ and ending in the
  interior of a non-horizontal edge.
\endroster

For each candidate system, one can construct a surface $F = F(\r_1,
\r_2, \r_3)$ in the exterior of $K$, called a {\it candidate surface}.
We refer the reader to [HO, P457] for the construction of the
surface.  Denote by $\hat F = \hat F(\r_1, \r_2, \r_3)$ the
corresponding closed surface obtained by capping off each boundary
component of $F$ with a disk.  Let $N(K)$ be a regular neighborhood of
$K$.  If $\delta$ is the boundary slope of $F$ on $\bdd N(K)$ then
$\hat F$ is a closed surface in the manifold $K_{\delta}$ obtained by
$\delta$ surgery on $K$.

When $\r_i$ ends at $\< \infty \>$ there may also be some ``augmented''
candidate surface, but fortunately this does not happen for Montesinos
knots of length 3.  The following is [HO, Proposition 1.1].  

\proclaim{Proposition 2.1}  Every incompressible,
$\partial$-incompressible surface in $S^3 - K$ having non-empty
boundary of finite slope is isotopic to one of the candidate surfaces.
\endproclaim

To find all toroidal surgeries on Montesinos knots of length 3, it
suffices to find all candidate systems $(\r_1, \r_2, \r_3)$ such that
$\hat F = \hat F(\r_1, \r_2, \r_3)$ is a torus.  By Theorem 1.2 all
toroidal $\hat F$ are incompressible.

\subhead 2.3.  The $u$-coordinate of a point and the length of an edgepath 
\endsubhead 

Any rational point $(x,y)$ in the diagram represents some curve system
$(a,b,c)$ on a 4-punctured sphere as shown in Figure 2.2.  (When $c$
is negative, reverse the tangency of the train track, and relabel $c$
by $-c$.)  The parameters $(a,b,c)$ and $(x,y)$ are related as
follows.
$$
\align
 y &= \frac c{a+b}  \\
 x &= \frac b{a+b}
\endalign
$$
See [HO, P455].  

\bigskip
\leavevmode

\centerline{\epsfbox{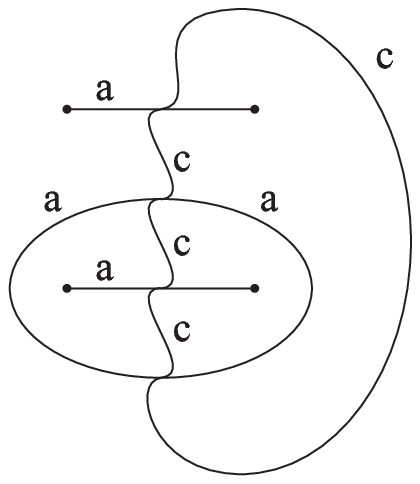}}
\bigskip
\centerline{Figure 2.2}
\bigskip

Note that $(a,b,c)$ is determined by $(x,y)$ up to scalar
multiplication, i.e., $(a,b,c)$ and $k(a,b,c)$ correspond to the same
rational point in $\cd$, so for any rational point $(x,y)$ one can
choose $a,b,c$ to be integers with $a > 0$.

A rational point in the
interior of an edge $\< p/q,\; r/s \>$ in $\cd$ corresponds to a curve
system $(1, b, c)$, which can be written as a linear
combination
$$
(1,b,c) = \a (1, s-1, r) + \b (1, q-1, p) 
$$
where $\a, \b$ are
positive rational numbers, and $\a + \b = 1$.  We write
$$
v =\a\< r/s \> + \b \< p/q \>
$$ 
to indicate that the point $v$ is related to $\< p/q \>$ and $\< r/s
\>$ as above.  The number $\a$ (resp.\ $\b$) is called the {\it
length\/} of the edge segment from $\< p/q \>$ (resp.\ $\< r/s \>$) to
$v$.  It is important to notice that this is not the euclidean length
of the segments of the edge cut by $v$, even if the length of the edge
is normalized to 1.  From the construction of the candidate surface
([HO, P455]), we see that traveling from the vertex $\< r/s
\>$ to the point $v$ above corresponds to adding $m\b$ saddles to the
surface, where $m$ is the number of times the surface intersects a
meridian of $K$, which must be an integer.  This fact will be useful
in the calculation of the Euler number of the resulting surface.

To make calculation easier, we introduce the {\it $u$-coordinate\/} of
a point $v$.  Define
$$ u = u(v) = \frac 1{1-x}$$ where $x$ is the $x$-coordinate of the
point $v$ in $\cd$.  Thus we have $x = (u-1)/u$.  The $u$-coordinate
has two important properties.  

(1) The $u$-coordinate of a vertex $\< p/q \>$ is $q$.

(2) The length of an edge segment is equal to its length in
$u$-coordinate when the length of the edge is normalized to 1, as
shown in the following lemma.

\proclaim{Lemma 2.2}
Let $v = \a\< r/s \> + \b\< p/q \>$.  Let $u=u(v)$, $u_0  = q$ and
$u_1 = s$ be
the $u$-coordinates of $v$, $\< p/q \>$ and $\< r/s \>$ respectively.  Then 
$$ u = \a u_1 + \b u_0.$$
In particular, $\a$ and $\b$ can be calculated by the following formulas.
$$ \align
\a & = \frac {u - u_0} {u_1 - u_0} = \frac {u-q}{s-q} \\
\b & = \frac {u_1 - u} {u_1 - u_0} = \frac {s-u}{s-q} 
\endalign
$$
\endproclaim

\proof
Suppose $\< p/q \>$ is represented by the curve system
$(a_0,b_0,c_0)$, and $\< r/s \>$ by $(a_1, b_1, c_1)$.  By definition
we may choose $(a_1, b_1, c_1)$ so that $a_1=a_0$.  Then the
$x$-coordinates of these points are $x_i = b_i/(a_i+b_i)$, hence $u_i
= 1/(1-x) = (a_i + b_i)/a_i = 1 + (b_i/a_i)$ for $i=0,1$.

By definition $v = a\< r/s \> + \b\< p/q \>$ is represented by $\< \a
a_1 + \b a_0, \a b_1 + \b b_0, \a c_1 + \b c_0\>$.  Using the facts that
$\a+\b=1$ and $a_0 = a_1$, we can calculate the $u$-coordinate of $\a
 \<r/s\> + \b \< p/q \>$ as follows.
$$\align
u & = 1 + \frac {\a b_1 + \b b_0}{\a a_1 + \b a_0} = 1 + \frac {\a b_1 + \b
b_0}{a_0} \\
  & = \a (1+\frac {b_1}{a_1}) + \b (1 + \frac {b_0}{a_0})  \\
  & = \a u_1 + \b u_0
\endalign
$$
\endproof

\definition{Definition 2.3}  The length $|\r|$ of an edgepath $\r$ in
$\cd$ is defined by counting the length of a full edge as 1, and
the length of a partial edge from $\< r/s \>$ to $\a\< r/s \> +
\b\< p/q \>$ as $\b$.
\enddefinition

\subhead 2.4. Euler numbers of points, edgepaths and surfaces 
\endsubhead

The knot $K = K(t_1, \; t_2, \; t_3)$ in $S^3$ can be constructed as
follows.  Let $(B_i, T_i)$ be a rational tangle of slope $t_i =
p_i/q_i$, where $B_i = D^2\times I$, and $T_i$ consists of two strings
with endpoints on the vertical diameters of $D^2 \times \bdd I$.
Gluing the end disks $D^2 \times \bdd I$ of the tangles in a cyclic
way, we get a knot in a solid torus $V$, which can be trivially
embedded in $S^3$ to produce the knot $K = K(t_1, \; t_2, \; t_3)$ in
$S^3$.  

Denote by $M_i$ the tangle space $B_i - \Int N(T_i)$.  Let $E_i$ be a
disk in $M_i$ separating the two arcs of $T_i$, and let $D_i = D_i^1
\cup D_i^2$ be a pair of disks properly embedded in $M_i$ such that
$D_i^j$ intersects the meridian of the $j$-th string of $T_i$ at a
single point and is disjoint from the meridians of the other string.

Define a number $m_i$ as follows.  If $\r_i$ is a not a constant path,
let $m_i$ be the minimal positive integer such that $m_i \times
|\r_i|$ is an integer.  If $\r_i$ is a constant path on $L(p_i/q_i)$
at a point with $u$-coordinate $\bu$, let $m_i$ be the smallest
positive integer such that $m_i \bu/q_i$ is an integer.  Let $n =
\lcm(m_1, m_2, m_3)$ be the least common multiple of $m_1, m_2, m_3$,
and let $m$ be a multiple of $n$.

\proclaim{Lemma 2.4}  Let $F(\r_i)$ be the surface in the tangle space
$M_i$ corresponding to the edgepath $\r_i$ constructed in [HO, P455]
which intersects a meridian of the tangle strands at $m$ points.

(1) If $\r_i$ is not a constant path then $\chi(F(\r_i)) = m (2 -
    |\r_i|)$.  

(2) If $\r_i$ is a constant path on the horizontal edge $L(t_i)$ with
    $u$-coordinate $\bu$, then $\chi(F(\r_i)) = m(1 + \bu/q_i)$

(3) Let $m_i$ be defined as above, and let $n = \lcm(m_1, m_2, m_3)$.
    If $F=F(\r_1, \r_2, \r_3)$ is connected and orientable, then it
    intersects a meridian of $K$ at either $n$ or $2n$ points.
\endproclaim

\proof (1) If $\r_i$ is not a constant path in the interior of
$L(t_i)$ then according to [HO, P457], $F(\r_i)$ is obtained from
$m$ copies of $D_i$ by adding some saddles.  For each full edge in
$\r_i$ one adds $m$ saddles, and for a partial edge of length $\b_i$
one adds $m\b_i$ saddles.  (By the choice of $m$ in the construction,
$m\b_i$ must be an integer.)  Since $\chi(mD_i) = 2m$, and adding a
saddle reduces the Euler characteristic by $1$, we have $\chi(F(\r_i))
= 2m - m|\r_i|$.

(2) Suppose $\r_i$ is a constant path in the interior of the
horizontal edge $L(t_i)$.  Then by [HO, P455], $F(\r_i)$ consists
of $m$ copies of $D_i$ and $k$ copies of $E_i$ for some $k$, hence
$\chi(F(\r_i)) = 2m + k$.  We need to determine the number $k$.

Let $(a', b', c')$ be the parameters of $\bdd D_i$ on the 4-punctured
sphere $\bdd M_i$.  Since it is a vertex on $L(t_i)$, we have $y_i
= p_i/q_i = c'/(a'+b')$, and $x_i = (q_i-1)/q_i = b'/(a'+b')$.  Also
$a' = 1$ because a meridian intersects $\bdd D_i$ at a single point.
Solving these equations gives $b' = q_i - 1$, and $c' = p_i$.

Let $(a'', b'', c'')$ be the parameters of $\bdd E_i$.  Then $a'' = 0$
because $\bdd E_i$ is disjoint from the meridians of $T_i$.  Examining
the curve on $\bdd B_i$ explicitly we see that $b'' = q_i$ and $c'' =
p_i$, hence it has parameters $(0, q_i, p_i)$.  

Now the parameters of $mD_i + kE_i$ are given by
$$m(1,\; q_i-1,\; p_i) + k(0,\; q_i,\; p_i) = (m,\; (m+k)q_i - m,\; (m+k)p_i).$$
Hence the $x$-coordinate and $u$-coordinate of the constant path
$\r_i$ satisfy
$$ \align
 x_i & = \frac {(m+k)q_i - m}{(m+k)q_i} \\
 \bu & = \frac 1{1-x_i} = \frac {(m+k)q_i}{m} = q_i + \frac {kq_i}{m}
\endalign
$$
Solving the last equation gives $k = m\bu/q_i - m$, therefore
$\chi(F(\r_i)) = 2m + k = m + m\bu/q_i$.

(3) By the proof of (2) the number $m$ must be a multiple of each
$m_i$, hence $m = kn$.  Then the number of initial disks and the
number of saddles in each step of the construction of $F(\r_1, \r_2,
\r_3)$ are all multiples of $k$.  One can therefore divide the initial
disks and the saddles in groups of $k$ sheets each.  When pinching
each group to a single sheet, we get a surface $F'$.  Thus $F$ lies in
a regular neighborhood of $F'$, intersecting each $I$-fiber exactly
$k$-times.  Since $F$ is orientable and connected, we see that $k=1$
or $2$, and $k=2$ if and only if $F'$ is non-orientable.  Therefore $m
\leq 2n$.
\endproof

\definition{Definition 2.5}  Let $v$ be a rational point in $\cd$ with
$u = u(v)$ as its $u$-coordinate, and let $\gamma$ be an edgepath with $v$ as
its ending point.  

\roster
\item If $v$ is not on a horizontal line, define $e(v) = \frac 13
(4-u(v))$.

\item If $v$ is on a horizontal line $L(p/q)$, define $e(v) = \frac
13 + u(v)(\frac 1q - \frac 13)$.

\item For an allowable edgepath $\r$ with ending point $v$, define
$e(\r) = e(v) - |\r|$. 

\item Given a candidate edgepath system $(\r_1, \r_2, \r_3)$, define
$\bar e = \bar e(\r_1, \r_2, \r_3) = \sum e(\r_i)$.
\endroster

The numbers $e(v)$, $e(\r)$ and $\bar e$ are called the {\it Euler
number\/} of a point, an edgepath, and a candidate system,
respectively.
\enddefinition

\example{Example 2.6}
(a) If $v \in L(p/2)$ then $e(v) = \frac 13 + \frac 16 u > 0$.

(b) If $v \in L(p/3)$ then $e(v) = \frac 13 > 0$.

(c) If $v \in L(p/q)$ and $q\geq 4$ then $e(v) = \frac 13 + u(\frac 1q
- \frac 13) \leq \frac 13 + q(\frac 1q - \frac 13) = \frac 43 - \frac
q3 \leq 0$, and $e(v) = 0$ if and only if $q=u=4$.

(d) At a vertex $v = \< p/q \>$, $u(v) = q$, so $e(v)$ can be
rewritten as $e(v) = \frac 13(4 - q)$; in particular, $e(\<p/q\>) <
0$ for all $q > 4$.
\endexample

\proclaim{Lemma 2.7}
Let $\r$ be an edgepath with $|\r| < 1$.  Let $v$ be the ending point
of $\r$, and let $u=u(v)$ be the $u$-coordinate of $v$.  Then

(1) $e(\r)>0$ if and only if (i) $v$ is on the horizontal edge $L(p/q)$
with $q \leq 3$, (ii) $v$ is on $\< p/1,\; r/s \>$ for some $s \leq 3$, or
(iii) $v$ is on $\< p/2,\; r/3 \>$ and $u > 2.5$.

(2) $e(\r)=0$ if and only if (i) $v$ is on $\< p/1,\; r/4 \>$, or (ii) $v$
    is on $\< p/2,\; r/3 \>$ and $u = 2.5$.
\endproclaim

\proof
Since $|\r|<1$, $\r$ cannot contain a full edge, hence if $v$ is a
vertex then $\r$ is a constant path.  By definition $\r$ is also a
constant path if $v$ is in the interior of a horizontal edge.  Thus if
$v$ is on a horizontal line $L(p/q)$ then the result follows from the
calculations in Example 2.6(4) because $e(\r) = e(v)$.  

We now assume that $v \in \< p/q,\; r/s \>$, and $v \neq \< p/q \>,
\<r/s\>$.  Let $\bu$ be the $u$-coordinate of $v$.  Then by definition
we have $e(\r) = \frac 13(4-\bu) - (s-\bu)/(s-q)$, which is a
linear function of $\bu$, $e(\r) = e(\< p/q \>)-1$ when $\bu=q$, and
$e(\r) = e(\< r/s \>)$ when $\bu=s$.  Thus when $s\geq 5$ we have $e(\r)
\leq 0$ at $\bu=q$, and $<0$ at $\bu=s$, hence $e(\r)<0$ for all $q<\bu\leq
s$.  The cases where $1\leq q < s\leq 4$ can be done one by one.  We omit
the details.
\endproof

The set of $v$ such that $e(\r)>0$ for some $\r$ ending at $v$ is
shown in Figure 2.3 for the square $[0,1] \times [0,1]$.  Those in the 
other squares are vertical translations of this graph.

\bigskip
\leavevmode

\centerline{\epsfbox{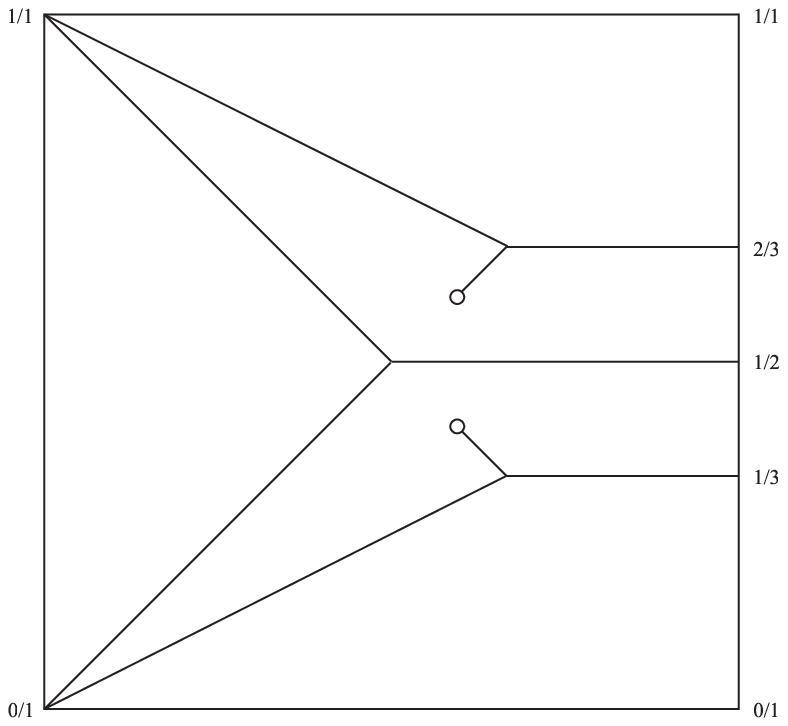}}
\bigskip
\centerline{Figure 2.3}
\bigskip

\proclaim{Theorem 2.8}  Let $K = K(p_1/q_1,\; p_2/q_2,\; p_3/q_3)$ be a
Montesinos knot of length 3, let $(\r_1, \r_2, \r_3)$ be a candidate
system, let $F=F(\r_1, \r_2, \r_3)$ be the associated candidate
surface, and let $\hat F = \hat F(\r_1, \r_2, \r_3)$ be the
corresponding closed surface.  Denote by $r = a/b$ the boundary slope of
$F$, where $a,b$ are coprime integers.  Then $\hat F$ is a torus if
and only if
$$ \bar e = \sum e(\r_i) = \frac {b-1}{b}
$$
In particular, if $r$ is an integer slope then $\bar e = 0$, and if
$r$ is a half integer slope then $\bar e = \frac 12$.
\endproclaim

\proof Let $(a_i,b_i,c_i)$ be the parameters of the ending point of
$\r_i$, chosen so that $a_1 = a_2 = a_3$ for all $i$, which will be
denoted by $m$.  Since $x_i = b_i/(b_i+a_i)$ are the same for a
candidate system, we have $b_1 = b_2 = b_3$, which we denote by $b$.

First consider the surface $F'$ obtained by gluing $F(\r_i)$ along the
three twice punctured disks $P_j$ on the boundary of the tangle
spaces.  Each $F(\r_i)$ intersects $P_j$ at $2m+b$ arcs, hence after
gluing the $F(\r_i)$ to each other, we have
$$\chi(F') = \sum \chi(F(\r_i)) - 3(2m + b).$$
By construction $F = F(\r_1, \r_2, \r_3)$ is obtained from
$F'$ by adding $2m+2b$ disjoint meridional disks in the solid torus
$S^3 - \cup B_i$, hence
$$\chi(F) = \chi(F') + 2m + 2b = \sum \chi(F(\r_i)) - 4m - b.$$
From $x = b_i/(a_i+b_i) = b/(m+b)$ and $\bu = 1/(1-x)$ one can solve
$b$ to obtain 
$$ b = \frac {mx}{1-x} = m (\bu-1).$$
When $r$ is an integer slope, we need to attach $m$ disks to $F$ to
obtain $\hat F$, hence 
$$\align
\chi(\hat F) & = \chi(F) + m = \sum \chi(F(\r_i)) - 3m - m(\bu-1) \\
& = m \sum (\frac 1m \chi(F(\r_i)) - \frac 23 - \frac 13 \bu)
\endalign
$$
If $\r_i$ is a constant path, by Lemma 2.4(2)
and Definition 2.5 we have
$$ \align
& \frac 1m \chi(F(\r_i)) - \frac 23 - \frac 13 \bu 
=  (1 + \frac {\bu}{q_i}) - \frac 23 - \frac 13 \bu  \\
= & \frac 13 + \bu (\frac 1{q_i} - \frac 13) =  e(\r_i) 
\endalign
$$
If $\r_i$ is not a constant path, by Lemma 2.4(1) and Definition 2.5 we
have 
$$ \align
& \frac 1m \chi(F(\r_i)) - \frac 23 - \frac 13 \bu 
= (2 - |\r_i|) - \frac 23 - \frac 13 \bu   \\
= & \frac 43 - \frac 13 \bu - |\r_i| =  e(\r_i)
\endalign
$$
Therefore we always have $\chi(\hat F) = m \sum e(\r_i)$, hence $\hat
F$ is a torus if and only if $\sum e(\r_i) = 0$.

The proof for $r=a/b$ and $b\neq 1$ is similar.  In this case $F$ has 
$m/b$ boundary components, so $\hat F$ is obtained by attaching $m/b$
disks to $F$.  Hence a similar calculation shows that
$$
\chi(\hat F) = \chi(F) + \frac mb = (\chi(F) + m) - \frac{b-1}{b}
m  = m \left( \sum e(\r_i) - \frac {b-1}{b} \right)
$$
Therefore in this case $\hat F$ is a torus if and only if $\sum
e(\r_i) = (b-1)/b$.
\endproof

\proclaim{Proposition 2.9}  Let $(\r_1, \r_2, \r_3)$ be a candidate system
such that $\hat F(\r_1, \r_2, \r_3)$ is a torus.  Suppose $\bu \leq
1$.  Then (i) $\bu = 1$, and (ii) $K = K(1/q_1,\; 1/q_2,\; 1/q_3)$ for
some (possibly negative) integers $q_i$, such that $|q_i|>1$, and at
most one $q_i$ is even.

The knots and the corresponding toroidal slopes are the same as those
in Theorem 1.1(1) and (2).  \endproclaim

\proof
Let $\r'_i$ be the part of $\r_i$ in the strip of $x \in [0,1)$ (i.e.\
$u \geq 1$), and let $y'_i$ be the ending points of $\r'_i$.  Then 
$|\r'_i|$
are all non zero integers, hence we have
$$ 0 \leq \bar e = (4-\bu) - \sum |\r_i| \leq 3 - \sum |\r'_i| \leq
0.$$ Thus all the inequalities above are equalities, and we have $\bu
= |\r_i| = |\r'_i| = 1$ for all $i$, so $\r_i = \r'_i$ contains only
one edge.  Since $\bu = 1$, $y_i = y'_i$ are integers.  By definition
of candidate system we have $\sum y_i = 0$, hence by choosing the
parameters properly we may assume that $y_i=0$ for all $i$.  It is now
easy to see that $K = K(1/q_1,\; 1/q_2,\; 1/q_3)$ for some $q_i$.
Since $K$ is of length 3, $|q_i| > 1$.  Since $K$ is a knot, at most
one $q_i$ is even.

The knots are the same as those in Theorem 1.1(1) and (2).  The
toroidal surface corresponding to a candidate system above is the
pretzel surface $S$, or its double cover if $S$ is nonorientable.  One
can draw the pretzel surface and show that the boundary slope of $F$
is the same as that in Theorem 1.1(1) and (2).  \endproof

Up to equivalence we may change the parameters of $K =
K(t_1,\; t_2,\; t_3)$ by the following moves.

(1) Replace all $t_i$ by $-t_i$;

(2) Permute $t_i$; 

(3) Replace $(t_1, t_2, t_3)$ by $(t_1+k_1, t_2+k_2, t_3+k_3)$, where
    $k_i$ are integers, and $\sum k_i = 0$.

If $(\r_1, \r_2, \r_3)$ is a candidate system for $K(t_1,\; t_2,\; t_3)$,
and $(t'_1, t'_2, t'_3)$ is equivalent to $(t_1, t_2, t_3)$ by the
above relations, then we can obtain a candidate system $(\r'_1, \r'_2,
\r'_3)$ for $K(t'_1,\; t'_2,\; t'_3)$ in the obvious way.  For example,
when $(t_1, t_2, t_3)$ is replaced by $(t_1+1, t_2-1, t_3)$, the
edgepath $\r'_1$ is obtained by moving $\r_1$ upward by one unit, and
$\r'_2$ downward by one unit.  Clearly the surface $F(\r_1, \r_2,
\r_3)$ is homeomorphic to $F(\r'_1, \r'_2, \r'_3)$.

Let $G = \< 0,\; -\frac 13 \> \cup \< 0,\; -\frac 12 \> \cup
\< -\frac 12,\; -\frac 13 \> \cup \<-\frac 12, -\frac 23 \> 
\cup \<-1,\; -\frac 12 \> \cup L(-\frac 12) \cup L(-\frac 13)$, as
shown in Figure 2.4. 

\bigskip
\leavevmode

\centerline{\epsfbox{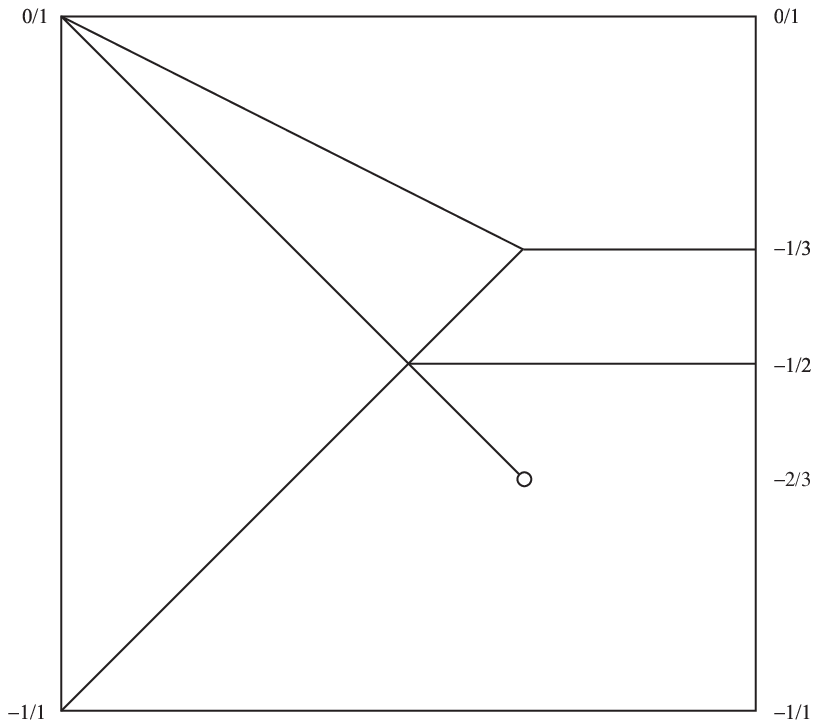}}
\bigskip
\centerline{Figure 2.4}
\bigskip

\proclaim{Lemma 2.10} Let $(\r_1, \r_2, \r_3)$ be a candidate edgepath
system for $K = K(t_1, t_2, t_3)$ such that the
corresponding surface $\hat F(\r_1, \r_2, \r_3)$ is a torus.  Let $v_i
= (x_i, y_i)$ be the ending point of $\r_i$, and let $\bu$ be the
$u$-coordinate of $v_i$.  Assume $\bu > 1$.  Then the following hold
up to re-choosing the parameters of $K$.

(1) $\sum y_i = 0$;

(2) $|y_i| + |y_j| \leq 1$ for any $i\neq j$;

(3) $0 < |y_i| \leq \frac 23$;

(4) $v_1$ is on the subgraph $G$ in Figure 2.4.
\endproclaim

\proof We may assume that the parameters of the knot has been chosen,
among equivalent knots, so that $\sum |y_i|$ is minimal.  The minimum
can be reached because (i) $\sum |y_i|$ remains the same when
permuting the parameters or replacing $(t_1, t_2, t_3)$ by $(-t_1,
-t_2, -t_3)$, and (ii) $\sum |y_i|$ goes to $\infty$ when $(t_1, t_2,
t_3)$ is replaced by $(t_1+k_1, t_2+k_2, t_3-k_1-k_2)$ and at least
one $k_i$ goes to $\infty$.

(1) This follows from the definition of candidate system.

(2) By permuting the $t_i$ and simultaneously changing their
signs if necessary, we may assume without loss of generality that
$-y_1 \geq y_2
\geq y_3 \geq 0$.  If the result is false then $-y_1 + y_2 > 1$.  But
then replacing $(t_1, t_2, t_3)$ of $K$ by $(t_1+1,
t_2-1, t_3)$ will give a  candidate system such that the
$y$-coordinates of the ending points are $y_1' = -y_1+1$, $y_2' =
y_2-1$, and $y_3' = y_3$, respectively.  One can check that $\sum
|y'_i| < \sum |y_i|$ if $|y_1|+|y_2|>1$.

(3) Since $\bu > 1$, $v_i$ cannot be on $L(0)$ as otherwise $\r_i$
would be a constant path on $L(0)$, so $t_i$ would be $0$,
contradicting the assumption that the parameters of $K = K(t_1, t_2,
t_3)$ are non-integers.  Therefore $y_i \neq 0$.  If $|y_1|> \frac
23$, say, then since $y_1 = -y_2 -y_3$, we would have $|y_i|>\frac 13$
for $i=2$ or $3$, which implies $|y_1|+|y_i|>1$, contradicting (2).

(4) Up to relabeling we may assume that $e(\r_1) \geq e(\r_i)$ for
$i=2,3$, and by taking the mirror image of $K$ if necessary we may
assume that $y_1<0$.  By Theorem 2.8,  $\bar e = \sum e(\r_i) \geq 0$,
hence either $e(\r_1)>0$, or $e(\r_i)=0$ for all $i$.

First assume $e(\r_1) > 0$.  Since $-1< y_1 <0$, by Lemma 2.7 $v_1$ is
on one of the edges in $G$, except that it may also be on the edge
$L(-\frac 23)$.  However, if $v_1 \in L(-\frac 23)$ then since $-y_1 =
y_2 + y_3$ and $-y_1 + y_i \leq 1$, we must have $y_2 = y_3 = \frac
13$, hence replacing $(t_1, t_2, t_3)$ by $(-t_2, -t_3, -t_1)$ will
give a new candidate system such that the ending point of the first
edgepath is on $L(-\frac 13)$, as required.

Now assume $e(\r_i) = 0$ for all $i$.  We may assume that no $v_i$ is
on $G$ or its reflection along the line $y=0$, as otherwise we may
choose the parameters of $K$ so that $v_1 \in G$.  Thus by Lemmas 2.7
and 2.10(3), each $v_i$ must be on $\<0, \pm \frac 14\>$.  However,
in this case one can show that $\sum y_i \neq 0$, contradicting Lemma
2.10(1).  Therefore this case cannot happen.  \endproof

\subhead 2.5.  Calculation of boundary slopes
\endsubhead

Denote by $e_-$ (resp.\ $e_+$) the number of edges in all the $\r_i$
on which a point moves downward (resp.\ upward) when traveling from
right to left.  Then the {\it twist number\/} of the edgepath system
$(\r_1, \r_2, \r_3)$ is defined as
$$
\tau = \tau(\r_1, \r_2, \r_3) = 2(e_- - e_+)
$$
Denote by $\delta = \delta (\r_1, \r_2, \r_3)$ the boundary slope of the
surface $F(\r_1, \r_2, \r_3)$.  The following lemma is due to Hatcher
and Oertel [HO], and can be used to calculate the boundary slope
$\delta$ for a given edgepath system.

\proclaim{Lemma 2.11}  Let $t_i$ be rational numbers, and let 
$(\r_1, \r_2, \r_3)$ be a candidate system with $\r_i$ starting at a
point on $L(t_i)$.  Then $\delta - \tau$ depends only on $t_i$ and is
independent of the paths $\r_i$.

Thus, if $F' = F(\r'_1, \r'_2, \r'_3)$ has boundary slope $\delta'$
and $\r'_i$ has starting point on $L(t_i)$, then $\delta= \tau +
\delta' - \tau'$, where $\delta' = \delta(\r'_1, \r'_2, \r'_3)$ and $\tau' =
\tau(\r'_1, \r'_2, \r'_3)$.  In particular if $F'$ is a Seifert
surface then $\delta = \tau - \tau'$.
\endproclaim

\proof 
This is on page Page 460 of [HO], where it was shown that $\delta =
\tau - \tau'$, where $\tau'$ is the twist number of the edgepath
system corresponding to the Seifert surface of $K$, starting from the
vertices $\<t_i\>$.  Therefore $\delta - \tau$ depends only on $(t_1,
t_2, t_3)$.
\endproof

\subhead 2.6. Notations and Conventions \endsubhead 

Throughout this paper we will denote by $\r_i$ the edgepath for the
$i$-th tangle, by $v_i$ the ending point of $\r_i$, by $\bar u$ the
$u$-coordinate of $v_i$, which must be the same for all $i$, and by
$y_i$ the $y$-coordinate of $v_i$.  Let $L$ be the union of the two
horizontal edges in $G$, i.e., $L = L(-1/2) \cup L(-1/3)$.

The case $\bu \leq 1$ has been discussed in Proposition 2.9.  Hence in
Sections 3--6 we will assume that $\bu > 1$.  By Lemma 2.10, in this
case we may choose the parameters $t_i$ of $K=K(t_1, t_2, t_3)$ to
satisfy the conclusions of that lemma; in particular, the ending point
$v_1$ of $\r_1$ lies on the subgraph $G$ of $\cd$ in Figure 2.4.  In
Sections 3--6 we will determine $K$ case by case, according to the
position of $v_1$ in $G$.

\head  3.  The case that $v_1 \in L$ and $\a_i = 0$ for $i=2$ or $3$  
\endhead

In this section we will discuss the case that one of the vertices, say
$v_1$, lies on the horizontal lines $L = L(-1/2) \cup L(-1/3)$, and
$\a_2 = 0$.  Note that the second condition is equivalent to that
either $\r_2$ is a constant path, or $v_2$ is a vertex of $\cd$.

\proclaim{Lemma 3.1}  $\gamma_i$ cannot all be constant paths.
\endproclaim

\proof  Let $y_i = p_i/q_i$ be the $y$-coordinates of the ending
points $v_i$ of $\r_i$, where $p_i, q_i$ are coprime integers.  By
Lemma 2.10(1) we have
$$\sum y_i = \sum \frac {p_i}{q_i} = 0 \tag{1.1} $$
If all $\r_i$ are constant paths, $K = K(p_1/q_1,\; p_2/q_2,\;
p_3/q_3)$ and since $K$ is a knot, at most one of the $q_i$ is even.
If one of the $q_i$, say $q_1$, is even, then we have
$$ q_2 q_3 p_1 + q_1 (q_3 p_2 + q_2 p_3) = 0.  $$ Since the first term
is odd and the other two are even, this is impossible.
If all $q_i$ are odd, then equation (1.1) implies that
$$ q_2 q_3 p_1 + q_1 q_3 p_2 + q_1 q_2 p_3 \equiv p_1 + p_2 + p_3
\equiv 0  \qquad \text{mod $2$}  $$
which implies that either one or three $p_i$ are even.  However, in
this case $K$ is a link of two components, which is again a contradiction.
\endproof

\proclaim{Lemma 3.2}  If $v_1$ is in the interior of $L$, and $v_2$ 
is in the interior of $L(p_2/q_2)$ for some $q_2 \leq 3$, then
$K =K(-1/2,\; 1/3,\; 1/(6+1/n))$ for some $n \neq 0, -1$, and $\bu =
6$.  
\endproclaim

\proof
Recall that we have assumed that $y_i$ satisfy the conclusions of
Lemma 2.10, hence $y_1 + y_2 +
y_3 = 0$, $0 < |y_i| \leq \frac 23$, and $|y_i| + |y_j| \leq 1$ for
$i\neq j$.

First assume $y_1 = -\frac 12$.  Then the above and the assumption of
$v_2 \in L(p_2/q_2)$ for $q_2 \leq 3$ imply that $y_2 = \frac 13$, and
$y_3 = -y_1 - y_2 = 1/6$.  The horizontal line $y= 1/6$ intersects the
graph $\cd$ at the horizontal edge $L(1/6)$ and one point on each edge
$\< 0,\; 1/q \>$ with $q\leq 6$.  By Lemma 3.1, $v_3$ cannot be in the
interior of $L(1/6)$ as otherwise we would have three constant paths.
It follows that $v_3$ must be on some $\< 0,\; 1/q \>$ with $q \leq
6$.  By calculating the intersection point of $y= 1/6$ with $\< 0,\;
1/q \>$ we see that $u \leq 3$ when $q \leq 5$, which would be a
contradiction because $v_2 \in \Int L(\frac 13)$ implies that $\bu>3$.
Therefore we must have $q = 6$, in which case $v_3$ is the vertex $\<
1/6 \>$.

By Definition 2.5  we have
$$\sum e(v_i) = (\frac 13 + 6(\frac 12 - \frac 13)) + \frac 13 + \frac
13(4-6) = 1$$
therefore by Theorem 2.8 we must have $\bar e = \sum e(\r_i) = \sum
e(v_i) - \sum |\r_i| = 0$, so there is exactly one edge in $\cup
\r_i$, which must be in $\r_3$ because $\r_1$ and $\r_2$ are constant
paths.  Therefore
$$K = K(-\frac 12,\; \frac 13,\; \frac 1{6+\frac 1n})$$ Since $\r_3$ must
be an allowable edgepath, we have $n\neq 0, -1$, and the result
follows.

Now assume $y_1 = -\frac 13$.  The case of $y_2 = \frac 12$ is similar
to the above, and we obtain the same knot up to equivalence.  If $y_2
= -\frac 12$ then $|y_2|+|y_3|>1$, contradicting Lemma 2.10.  In all
other cases we have $v_i \in \Int L(p_i/3)$ for each $i$, which implies
that $\r_i$ are all constant paths, contradicting Lemma 3.1.
\endproof

\proclaim{Lemma 3.3}  Suppose $v_1 \in L(-1/3)$.  Then
$v_2$ cannot be in the interior of a horizontal edge $L(p_2/q_2)$ with
$q_2\geq 4$.
\endproclaim

\proof
If this is not true then $\r_2$ is a constant path, so by Lemma 3.1
$\r_3$ cannot be a constant path, hence $|\r_3| > 0$.  By Definition
2.5 we have
$$\align
e(\r_1) & = \frac 13 \\
e(\r_2) & = \frac 13 + \bu (\frac 1{q_2} - \frac 13) \\
e(\r_3) & = \frac 13 (4-\bu) - |\r_3| \\
0 \leq \bar e & = \sum e(\r_i) = 2 - \frac 5{12} \bu - \bu(\frac 14 -
\frac 1{q_2}) - |\r_3|  \\ 
& \leq  2 - \frac 5{12} \bu - |\r_3| < 2 - \frac 5{12} \bu 
\endalign
$$
This gives $5 > \bu$.  Since $\bu \geq q_3$ and $q_3 \geq 4$, we must
have $q_2 = 4$.  Assume
$v_3$ is on the edge $\< p_3/q_3,\; t_3/s_3 \>$.  Then $s_3 > \bu \geq 4$,
so $s_3 \geq 5$.

Define $\b_3(u) = (s_3 - u)/(s_3 - q_3)$.  Then $\b_3(\bu)$ is the
length of the last edge segment in $\r_3$, so $\b_3(\bu) \leq |\r_3|$. 

The function $e(u) = 2 - \frac 5{12} u - \b_3(u)$ is a linear function
of $u$.  We have $e(5) < 0$, and $e(\bu) \geq \bar e \geq 0$, for some
$4<\bu<5$, so $e(4)>0$, and hence $\b_3(4) < \frac 13$.  Since
$\b_3(4) = (s_3-4)/(s_3-q_3)$, this is true if and only if $q_3 = 1$
and $s_3 = 5$.  Hence $v_3$ is on an edge $E_3 = \< p_3/1,\; t_3/5 \>$
for some $p_3, t_3$.  Since $0<|y_3|
\leq \frac 23$, we must have $E_3 = \<0, \; \pm 1/5 \>$.

By assumption we have $y_1 = -\frac 13$.  Since $|y_2|=|p_2/q_2| =
|p_2|/4$ and $|y_1|+|y_2|\leq 1$, we must have $|y_2| = \frac 14$,
hence $(y_2, y_3) = (-\frac 14, \frac 5{12})$ or $(\frac 14, \frac
1{12})$.  It is easy to see that the horizontal line $y= \frac 5{12}$
does not intersect the edge $E_3$ above.  Therefore we must have
$(y_2, y_3) = (\frac 14,\; \frac 1{12})$, hence $E_3 = \<0, \frac
15\>$.  

The line equation of $E_3$ is given by $y =  x/4$, hence the only
solution for $y_3 = 1/12$ and $x>0$ is at $x = 1/3$, which has 
$u$-coordinate $u = 1/(1-x) = 3/2 <4$.  Therefore there is no
solution in this case because $\bu > 4$.
\endproof

\proclaim{Lemma 3.4}  Suppose $v_1 \in L(-\frac 12)$.  Then $v_2$
cannot be in the interior of a horizontal edge $L(p_2/q_2)$ with
$q_2\geq 4$.
\endproclaim

\proof
Similar to Lemma 3.3, we have
$$\align
e(\r_1) & = \frac 13 + \bu(\frac 12- \frac 13) = \frac 13 + \frac 16 \bu \\
e(\r_2) & = \frac 13 + \bu (\frac 1{q_2} - \frac 13) \\
e(\r_3) & = \frac 13 (4-\bu) - |\r_3|  \\
0 \leq \bar e & = \sum e(\r_i) = 2 - \bu(\frac 12 - \frac 1{q_2}) - |\r_3|
\endalign
$$

Since $K$ is a knot, $q_2$ must be odd, hence $q_2 \geq 5$.  Hence
from the above we have $2 - \bu(\frac 12 - \frac 15) \geq 0$, so $\bu <
7$.  Since $q_2 < \bu$, we must have $q_2 = 5$.  Therefore $y_2 =
p_2/q_2 = \pm k/5$ for some $k=1,2,3,4$.  Since $|y_i| + |y_j| \leq 1$ and
$|y_1| = \frac 12$, we must have $y_2 = \frac 15$ or $\frac 25$.  (The
cases of $y_2 = -\frac 15$ and $-\frac 25$ are impossible because then
$y_3 > \frac 12$, so $|y_1| + |y_3| > 1$.)

When $y_2 = \frac 25$, we have $y_3 = -y_1 -y_2 = \frac 1{10}$.  The
intersection of the line $y = \frac 1{10}$ and $\cd$ is the union of
$L(\frac 1{10})$ and one point in each edge $\< 0,\; 1/q_3 \>$ for $q_3
\leq 10$.  Since $\bu < 7$, $v_2$ cannot be on $L(\frac 1{10})$.  
By direct calculation we see that the $u$ value of the intersection
between $y=\frac 1{10}$ and $\< 0,\; 1/q_3 \>$ is
$$
u = \frac 1{1-(q_3 -1)y} = \frac {10}{11 - q_3}
$$
which gives $u\leq 5$ for $q_3 \leq 9$, and $u = 10$ for $q_3 = 10$.
Since $5<\bu<7$, there is no solution in this case.

When $y_2 = \frac 15$, we have $y_3 = -y_1 - y_2 = \frac 3{10}$.  The
horizontal line $y=\frac 3{10}$ intersects $\cd$ at $L(\frac 3{10})$,
and one point on each of $\< 0,\; \frac 12 \>$, $\< 0,\; \frac 13 \>$, $\<
\frac 13,\; \frac 14 \>$, and 
$\< \frac 13,\; \frac 27 \>$.  (These are all the edges $\< t_1,\; t_2
\>$ with $t_1$ and $t_2$ on opposite sides of $y=\frac 3{10}$.)  As
above, one can calculate the $u$-coordinate of the intersection to
show that there is no intersection point on the interval $5<u<7$.
Hence there is no solution in this case either.
\endproof

We now assume that $v_1 \in L$, and $v_2, v_3$ are not in the interior
of horizontal edges.  Since $\a_2 = 0$, the ending point $v_2$ of
$\gamma_2$ must be a vertex of $\cd$.  The following two lemmas determine
all knots with this property.

\proclaim{Lemma 3.5} Suppose $v_1 \in L(-1/3)$, $v_2$ is a vertex
of $\cd$, and $v_3$ is not in the interior of a horizontal edge.  Then
$K = K(-1/3,\; -1/(3+1/n),\; 2/3)$ for some odd $n \neq -1$, and $\bu
= 3$.
\endproclaim

\proof
By Definition 2.5 we have 
$$0 \leq \bar e = \sum e(\r_i) \leq \frac 13 + 2 \times \frac 13
(4-\bu) - \sum |\r_i|,$$ which gives $\bu \leq 4.5$.  Since $v_2$ is
vertex and $v_1 \in L(-1/3)$, $\bu$ is an integer, and $\bu \geq 3$.
Hence $\bu = 3$ or $4$.  Thus $v_2 = \< \pm 1/3 \>$, $\< \pm 2/3 \>$,
$\< \pm 1/4 \>$ or $\< \pm 3/4 \>$.  Since $|y_i| \leq \frac 23$, we
cannot have $|y_2| = 3/4$.

When $y_2 = -1/3$, we have $y_3 = 2/3$, so all the three $v_i$ are
vertices, and $\sum e(v_i) = 1$.  Since $\bar e = \sum e(v_i) - \sum
|\r_i|$, by Theorem 2.8 we must have $\sum |\r_i| = 1$, so there is
one full edge in some $\r_i$.  Because of symmetry $K$ is equivalent
to $K(-1/3,\; -1/(3+1/n),\; 2/3)$.  Since $\r_3$ is allowable, $n\neq
-1$, and since $K$ is a knot, $n$ must be odd.  This gives the knots
listed in the lemma.

When $y_2 = 1/3$ or $-2/3$, $y_3 = 0$ or $1$, which is not a solution.
When $y_2 = 2/3$ we have $y_3 = -1/3$, which gives the same solution
as above.

When $y_2 = \frac 14$, we have $y_3 = -y_1 - y_2 = \frac 1{12}$.  Now
$v_3$ lies on the intersection of $y=\frac 1{12}$ and $u=4$, which
is a point on the edge $\< 0,\; 1/10 \>$.  We have 
$$ 0 \leq \bar 3 \leq \sum e(v_i) = \frac 13 + 0 + (-\frac 23) < 0$$
hence it is not a solution.

When $y_2 = -\frac 14$, $y_3 = \frac 13 + \frac 14 = \frac 7{12}$.
One can check that the point of $(u,y) = (4, \frac 7{12})$ lies on the
edge $\< \frac 12, \frac 35 \>$.  We have $\b_3 = (\bu - 2)/(5-2) =
1/3$.  Since $\b_1 = \b_2 = 0$, by Lemma 2.11 the boundary slope of the
surface is $\delta \equiv 2(e_- - e_+) \equiv \pm 2\b_3 = \pm 2/3$ mod
$1$, hence by Theorem 2.8 we have $\bar e = 2/3$.  On the other hand,
we have  
$$
\bar e = \sum e(\r_i) = \frac 13 + 0 + (\frac 13(4-4) - \frac 13) = 0$$
This comtradiction completes the proof of the lemma.
\endproof

\proclaim{Lemma 3.6} Suppose $v_1 \in L(-\frac 12)$, $v_2$ is a vertex of
$\cd$, and $v_3$ is not in the interior of a horizontal edge.  Then
$K$ and $\bu$ are given by one of the following.

(i) $K(-1/2,\; 1/5,\; 2/7)$, $\bu = 5$;

(ii) $K(-1/2,\; 2/5,\; 1/9)$, $\bu = 5$; 

(iii) $K(-1/2,\; 1/5,\; 1/(3+1/n))$, $n$ even, $n\neq 0$, $\bu = 3$;

(iv) $K(-1/2,\; 1/3,\; 1/(5+1/n))$, $n$ even, $n\neq 0$, $\bu = 3$.
\endproclaim

\proof  Let $v_2 = \<y_2\> = \<p_2/q_2\>$.  By definition of $e(\r_i)$
and Theorem 2.8 we have
$$
\align
0 \leq \bar e = \sum e(\r_i) & = (\frac 13 + \frac 16 \bu) + 2 \times
\frac 13 (4-\bu) - |\r_2| - |\r_3|   \\
& = 3 - \frac 12 \bu - |\r_2| -|\r_3|
\endalign
$$ 
which gives $ \bu = q_2 \leq 6$.  We have $y_2 > 0$ since otherwise
$y_3 = -y_1 - y_2 > \frac 12$, so $|y_1| + |y_3| > 1$, contradicting
Lemma 2.10(2).  Similarly we must have $y_2 \leq \frac 12$ as
otherwise we would have $|y_1| + |y_2| > 1$.  Moreover, $y_2 \neq
\frac 12$ as otherwise we would have $y_3 = -y_1 -y_2 = 0$,
contradicting Lemma 2.10(3).  Therefore $y_2 = \frac 13$, $\frac 14$,
$\frac 15$ $\frac 25$, or $\frac 16$.  In each case $v_3$ is uniquely
determined by the facts that $u(v_3) = \bu = q_2$, and $y_3 = -y_1 -
y_2$.  We separate the cases.

\medskip
CASE 1.  $y_2 = \frac 13$.  

We have $\bu = 3$, and $y_3 = \frac 16$.  The point $v_3$ lies on the
edge $\< 0, \frac 15 \>$, with $\b_3 = (5-3)/(5-1) = 1/2$, which is
the length of the last segment of $\r_3$.  Hence
$$\sum e(v_i) - \b_3 = \frac 56 + \frac 13 + \frac 13 - \frac 12 = 1$$
so there is an extra edge, whose ending point is either $\< 1/3 \>$ or
$\< 1/5 \>$.  Since $K = K(-\frac 12,\; a_2/b_2,\; a_3/b_3)$ is a
knot, the numbers $b_2$ and $b_3$ must be odd.  Combining these, we
see that $K$ is equivalent to a knot of type (iii) or (iv) in the
Lemma.

\medskip
CASE 2.  $y_2 = \frac 14$.

We have $\bu = 4$ and $y_3 = -y_1 - y_2 = \frac 14$, so $v_3$ is also
at $\< \frac 14 \>$.  Since $ \sum e(v_i) = (\frac 13 + \frac 16 \bu) + 0
+ 0 = 1$, there is one extra edge.  It follows that $K = K(-\frac 12,\;
\frac 14,\; 1/(4+\frac 1n))$, which is a link of at least two
components.  Therefore there is no solution in this case.

\medskip
CASE 3.  $y_2 = \frac 15$.  

We have $\bu = 5$ and $y_3 = \frac 12 - \frac 15 = \frac 3{10}$.  The
vertex $v_3$ lies on the edge $\< \frac 13, \frac 27 \>$, and $\b_3 =
(7-5)/(7-3) = 1/2$.  We have $\sum e(v_i) - \b_3 = (\frac 13 + \frac
56) + \frac 13(4-5) + \frac 13(4-5) - \frac 12 = 0$, so there is no
extra edge.  The knot is $K(-\frac 12,\; \frac 15,\; \frac 27)$.

\medskip
CASE 4.  $y_2 = \frac 25$.

Then $\bu = 5$ and $y_3 = \frac 1{10}$.  The point $v_3$ lies on the
edge $\< 0, \frac 19 \>$, and $\b_3 = (9-5)/(9-1) = 1/2$.  We have
$\sum e(v_i) - \b_3 = (\frac 13 + \frac 56) + \frac 13(4-5) + \frac
13(4-5) - \frac 12 = 0$, so there is no extra edge.  The knot is
$K(-\frac 12,\; \frac 25,\; \frac 19)$.

\medskip
CASE 5.  $y_2 = \frac 16$.

Then $\bu = 6$ and $y_3 = \frac 12 - \frac 16 = \frac 13$.  
The point $v_3$ is in the interior of the horizontal line $L(\frac
13)$, which contradicts the assumption.  Therefore there is no
solution in this case.
\endproof

\proclaim{Proposition 3.7}  Suppose $v_1 \in L$ and $\a_i = 0$ for
$i=2$ or $3$.  Then $K$ is equivalent to one of the knots listed in
Lemma 3.2, 3.5 or 3.6.  
\endproclaim

\proof  
By symmetry we may assume that $\a_2 = 0$, so $v_2$ is either a vertex
or in the interior of a horizontal line $L(p_2/q_2)$.  The first case
is covered by Lemmas 3.5 and 3.6.  In the second case by Lemmas 3.3
and 3.4 we must have $q_2 \leq 3$.  Since $y_2 \neq 0$, we have $q_2
\geq 2$.  We may now apply Lemma 3.2 unless $v_1$ is a vertex of $L$,
which happens only if $q_2 = 2$ and $v_1 = \<-\frac 13\>$.  If that is
the case, then we may consider the (equivalent) knot $K(-p_2/q_2,\;
-p_1/q_1,\; -p_3/q_3)$ instead, which has the property that the ending
points of the corresponding edgepaths are $(v'_1, v'_2, v'_3)$, with
$v'_1 \in L(-\frac 12)$ and $v'_2 = \< \frac 13 \>$, which has been
covered by Lemma 3.6.
\endproof

\head 4.  The case that $v_1 \in L(-\frac 13)$ and $\a_i \neq 0$ 
for $i=2,3$
\endhead 

In this section we will assume that $v_1 \in L(-\frac 13)$, and $v_i$
is in the interior of a non-horizontal edge $E_i = \< p_i/q_i,\;
r_i/s_i \>$ and hence $0 < \a_i < 1$ for $i=2,3$.  We have $\bu =
u(v_1) \geq 3$.

Define 
$$ \align
\b_i(u) & = \frac {s_i - u}{s_i - q_i}  \\
e(u) & = 3 - \frac 23 u - \b_2(u) - \b_3(u).
\endalign
$$
Then $e(u)$ is a linear function of $u$.  
Recall that  $\b_i = (s_i - \bu)/(s_i - q_i) = \b_i(\bu)$.  We have

\proclaim{Lemma 4.1} Suppose $v_1 \in L(-\frac 13)$, and $\a_i \neq 0$
for $i=2,3$.  Then

(1)  $0 \leq \bar e \leq  e(\bu)$, and 

(2)  $3 \leq \bu < 4.5$.
\endproclaim

\proof
(1) By definition we have 
$$ \align
0 \leq \bar e & = \sum e(\r_i) = \frac 13 + 2 \times (\frac 13(4-\bu)) -
|\r_2| - |\r_3| \\
& \leq  3 - \frac 23 \bu - \b_2 - \b_3 = e(\bu).
\endalign
$$

(2) Since $v_i$ is in the interior of non-horizontal edge, $\b_i > 0$
for $i=2,3$.  Hence the above inequality implies that $\bu < 4.5$.
Since $v_1 \in L(-\frac 13)$, we have $\bu \geq 3$.
\endproof

\proclaim{Lemma 4.2} Suppose $v_1 \in L(-\frac 13)$, and $\a_i \neq 0$
for $i=2,3$.  Then $q_i \leq 3$ for $i=2,3$. 
\endproclaim

\proof  Since $q_i \leq \bu$, by Lemma 4.1(2) we have $q_i \leq 4$ for
$i=2,3$.  If $q_2 = 4$ then by Lemma 4.1(1) we have $4 < \bu < 4.5$, so
$\b_2 = (s_2 - \bu)/(s_2 - 4) \geq 1/2$.   Therefore $e(\bu) = 3 -
\frac 23 \bu - \b_2 - \b_3 < 0$, which is a contradiction to Lemma
4.1(1).  
\endproof

\proclaim{Lemma 4.3} Suppose $v_1 \in L(-\frac 13)$ and $\a_i \neq 0$
for $i=2,3$.  Then $q_2 \leq 2$. 
\endproclaim

\proof
Assume to the contrary that $q_2 > 2$.  Then by Lemma 4.2 we must have
$q_2 = 3$.  In this case $\b_2(3) = -1$, so $e(3) = 3 - \frac 23 (3) -
1 - \b_3 < 0$.  

First assume $s_2 \geq 5$.  If $s_3 \geq 5$ then $\b_i(5) \geq 0$ for
$i=2,3$, so $e(5) = 3 - \frac 23 \times 5 - \b_2(5) - \b_3(5) <0$, and by
linearity we have $e(\bu) < 0$, which is a contradiction.  If $s_3 =
4$ then $3<\bu<4$.  Since $\b_2(4) \geq \frac 12$ we have $e(4) < 3 -
\frac 23 4 - \b_2(4) < 0$, which again contradicts the fact that
$e(\bu) \geq 0$.  

We may now assume $s_2 = 4$.  Then $3<\bu<4$.  We have $e(3) < 0$, and
$e(\bu)\geq 0$, hence $e(4) = \frac 13 - \b_3(4) > 0$, i.e.\ $\b_3(4)
= (s_3 - 4)/(s_3 - q_3) < \frac 13$.  

If $q_3 = 3$ then $\b_3(4) < \frac 13$ implies that $s_3 = 4$.  Since
$|y_i| \leq \frac 23$, $E_3 \neq \< \pm \frac 23, \pm \frac 34 \>$, so we
must have $E_3 = \< \pm \frac 13, \pm \frac 14 \>$.  Since $(q_2, s_2)
= (3,4)$, the same is true for $E_2$, hence $\frac 14 < |y_i| <
\frac 13$ for $i=2,3$.  Therefore either $y_2 + y_3 = 0$ (if $E_2 \neq
E_3$), or $|y_2 + y_3| > \frac 12$, either case contradicting the fact
that $y_2 + y_3 = -y_1 = \frac 13$.

If $q_3 = 2$ then since $s_3$ is coprime with $q_3$, $s_3 \neq 4$, so
we must have $s_3 \geq 5$.  Thus $\b_3(4) = (s_3 - 4)/(s_3 - 2) \geq
1/3$, which is a contradiction.  

If $q_3 = 1$ then $\b_3(4) = (s_3 - 4)/(s_3 - 1) < 1/3$ implies that
$s_3 = 4$ or $5$, so $E_3 = \< 0,\; \pm \frac 14 \>$ or $\< 0,\; \pm
\frac 15 \>$.  As above, $E_2 = \<\pm \frac 13, \pm \frac 14\>$, and
one can check that there is no solution to the equation $y_1 + y_2 +
y_3 = 0$ in these cases.
\endproof

\proclaim{Proposition 4.4}
Suppose $v_1 \in L(-\frac 13)$, and $\a_i \neq 0$ for $i=2,3$.  Then
$K$ is equivalent to either $K(-1/3,\; 1/5,\; 1/5)$ or $K(-1/3,\;
1/4,\; 1/7)$, and $\bu = 3$.
\endproclaim

\proof
First assume that $s_2, s_3 \geq 5$.  Then $\b_i(5) \geq 0$, hence
$e(5) = 3 - \frac 23 \times 5 - \b_2(5) - \b_3(5) < 0$.  By Lemma 4.1
we have 
$$ e(3) = 3 - \frac 23 \times 3 - \b_2(3) - \b_3(3) = 1-\b_2(3)
- \b_3(3)$$
Since $e(\bu) \geq 0$ for some $3\leq \bu < 5$, by linearity we have
$e(3) \geq 0$.  Therefore, one of the $\b_i$, say $\b_2$, satisfies 
$$\b_2(3)= \frac {s_2 - 3}{s_2 - q_2} \leq \frac 12$$ Since $s_2 \geq
5$, this is true if and only if $(q_2, s_2) = (1,5)$, in which case
$\b_2(3) = \frac 12$, hence $\b_3(3) \leq \frac 12$, and the same
argument as above shows that $(q_3, s_3) = (1,5)$.  It follows that $K
= K(-1/3,\; 1/5,\; 1/5)$, and $\bu =3$ because $e(\bu) < 0$ if $\bu >
3$.

Now assume $s_2 = 4$.  Then $3\leq \bu < 4$.  Since $q_2$ is coprime
with $s_2$ and $q_2 < 3$ by Lemma 4.3, we must have $q_2 = 1$, so 
$(q_2,s_2) = (1,4)$.  Since $|y_2|\leq \frac 23$, $E_2 \neq \<\pm 1,
\pm \frac 34\>$.  Therefore we must have $E_2 = \<0, \pm \frac 14\>$.  

Suppose $E_2 = \<0, -\frac 14\>$.  Then $- E_1 - E_2$ coincides with
the edge $\<\frac 12, \frac 35\>$ on $3\leq u\leq 4$, so $E_3 =
\<\frac 12, \frac 35\>$ if there is a solution.  By Example 2.6 and
Definition 2.5 we have
$$
\bar e = \sum e(\r_i) \leq \sum e(v_i) - \b_3 = 
\frac 13 + 0 + \frac 13(4-\bu) -
\frac {5-\bu}{5-2} = 0$$
On the other hand, since $E_2$ goes upward and $E_3$ downward when
traveling from right to left, we have $e_- \equiv \b_3$ and $e_+
\equiv \b_2$ mod $1$, hence by Lemma 2.11 the boundary slope of the
surface $\delta$ satisfies
$$\delta \equiv 2(e_- - e_+) \equiv 2 (\b_3 - \b_2)  
= 2 (\frac {5-\bu}{5-2} - \frac {4-\bu}{4-1}) = \frac 23 
$$
It follows from Theorem 2.8 that $\bar e = \frac 23$, which is a
contradiction.  Therefore there is no solution in this case.

Now assume $E_2 = \<0, \frac 14 \>$.  Then $y_2 = x/3$, and $y_3 =
-y_1 - y_2 = (1-x)/3$.  One can check that when $\bu = 3$, $v_3$ is on
the edge $\<0, \frac 17\>$, in which case we have $\bar e = 0$, and
$K = K(-\frac 13,\; \frac 14,\; \frac 17)$.  We need to show that
there is no solution when $\bu>3$.

The line segment $y = (1 - x)/3$ is below the line $y=1-x$, hence from
Figure 2.1 we see that $y_3$ is on an edge $E_3 = \<0,\; 1/s_3\>$ for
some $s_3$.  Since the line segment has negative slope, and since it
intersects $\<0, \frac 17\>$ at $u = 3$, we must have $s_3>7$ when $u>3$.
By definition we have
$$
\bar e = \sum e(\r_i) = \frac 13 + 0 + (\frac 13(4-\bu) - \frac {s_3 -
\bu}{s_3 - 1}) 
$$
For $s_3>7$, the right hand side is negative
for $\bu = 3$ and $4$, hence by linearity it is negative for all
$3\leq \bu \leq 4$.  By Theorem 2.8 there is no solution in
this case.
\endproof

\head 5.  The case that $v_1 \in L(-\frac 12)$ and $\a_i \neq 0$ 
for $i=2,3$  
\endhead 

In this section we will assume that $v_1 \in L(-\frac 12)$, and $v_i$
is in the interior of a non-horizontal edge $E_i = \< p_i/q_i,\;
r_i/s_i \>$ and hence $0 < \a_i < 1$ for $i=2,3$.  By Lemma 2.10 (2)
and (3) we must have $0<y_i < \frac 12$ for $i=2,3$, hence $0\leq z
\leq \frac 12$ when $z= p_i/q_i$ or $r_i/s_i$ and $i=2,3$.

As before, define $\bar e = \sum e(\r_i)$, and $\b_i(u) = (u -
q_i)/(s_i - q_i)$.  Let $\b_i = \b_i(\bu)$.  Define a function
$$ e(u) = 3 - \frac 12 u - \b_2(u) - \b_3(u).$$  Note that this is
different from the function $e(u)$ defined in Section 4.

Define $l_i = s_i - q_i$.  Given an edge $E$ in $\cd$ and a number
$t$, we use $t - E$ to denote the set of points $\{t-t' \; | \; t' \in
E\}$.

\proclaim{Lemma 5.1}  Suppose $v_1 \in L(-\frac 12)$ and $\a_i \neq 0$
for $i=2,3$.   Then

(1) $$
\align
0 \leq \bar e & \leq e(\bu) = 3 - \frac 12 \bu - \b_2 - \b_3 = 1 -
\frac 12 \bu + \a_2 + \a_3 \\ 
  & = 1 - \frac 12 \bu + \frac {\bu - q_2} {l_2} +  \frac {\bu - q_3} {l_3} \\
  & = (1 - \frac{q_2}{l_2}  - \frac{q_3}{l_3}) 
     - (\frac 12 - \frac{1}{l_2}  - \frac{1}{l_3}) \bu
\endalign
$$

(2)  $2 \leq \bu < 6$.
\endproclaim

\proof (1) By Theorem 2.8 we have $\bar e = \sum e(\r_i) \geq 0$.
Since $v_1 \in L(-\frac 12)$, by Definition 2.5 and Example 2.6 we
have 
$$\bar e = (\frac 13 + \frac 16 \bu) + 2 \times \frac 13(4-\bu) -
|\r_2| - |\r_3| \leq e(\bu)$$  
The other equalities are just different expressions of $e(\bu)$.

(2)  Since $v_1 \in L(-\frac 12)$, we have $\bu \geq 2$.  Since $v_i$
     is in the interior of non-horizontal edges, we have $\b_i > 0$,
     so $0 \leq e(\bu) < 3 - \frac 12 \bu$, hence $\bu < 6$.
\endproof

\proclaim{Lemma 5.2} Suppose $v_1 \in L(-\frac 12)$, and $\a_i\neq 0$
     for $i=2,3$.  Then $q_i \leq 3$ for $i=2,3$.
\endproclaim

\proof  Note that $q_i \leq \bu < 6$.  
If $q_2 = 5$ then $\bu \in [5, 6)$.  We have $e(6) \leq 3 - \frac 12
\times 6 = 0$, and
$e(5) = 3 - \frac 52 - \b_2(5) - \b_3(5) < 0$ because $\b_2(5) = 1$.
Since $e(u)$ is linear, $e(\bu) < 0$, a contradiction.  Therefore we
may assume that $q_2 = 4$.

First assume $s_2 > 5$.  We have $e(4) = 1 - \b_2(4) - \b_3(4)
= -\b_3(4) < 0$ and $e(\bu)
\geq 0$.  If $s_3 > 5$ then $\b_i(6) \geq 0$, so $e(6) \leq 3 - \frac
12 \times 6 = 0$, which contradicts $4<\bu < 6$ and the linearity of $e(u)$.
If $s_3 = 5$ then $4 < \bu < 5$, and $e(5) = 3 - \frac 12 \times 5 - \b_2(5)
\leq 0$, which again is a contradiction.

Now assume $s_2 = 5$.  Then $E_2 = \< p_2/4,\; r_2/5 \>$.  Since
$0<y_2<\frac 12$, we must have $E_2 = \< \frac 14,\; \frac 15 \>$.
However, one can check that the interior of the line segment $\frac 12
- E_2$ lies in the interior of the triangle with vertices $\<\frac
13\>$, $\<\frac 14 \>$, $\< \frac 27 \>$, hence there is no solution
to the equation $\sum y_i = 0$ for $v_2 \in \Int E_2$.  This
contradiction completes the proof of the lemma.
\endproof

\proclaim{Lemma 5.3} If $v_1 \in L(-\frac 12)$, $\a_i > 0$ for
$i=2,3$, and $q_3 \leq q_2 = 3$, then $K = K(-1/2,\; 2/5,\; 1/7)$ and
$\bu = 4$, or $K(-1/2,\; 1/5,\; 2/7)$ and $3<\bu < 5$.
\endproclaim

\proof
Since $y_2 \in (0, \frac 12)$, we must have $E_2 = \< 1/3,\;
r_2/s_2\>$.  From Figure 2.1 we see that $r_2/s_2 \geq \frac 14$, hence
$y_2 > \frac 14$, and $y_3 = - y_1 - y_2 < \frac 14$.  Since $q_3 \leq
3$ and $y_3 < 1/4$, from Figure 2.1 we see that $p_3/q_3 = 0$, so $E_3
= \< 0,\; 1/s_3\>$ for some $s_3 \geq 4$.  As before, put $l_i = s_i -
q_i$.

Since $y_1 = -\frac 12$ and $y_3 = -y_1 - y_2 = \frac 12 - y_2$, we
see that $v_3$ lies on the intersection of $E_3$ and $\frac 12 - E_2$.
When $E_2 = \< 1/3,\; 1/4\>$, one can check that the interior of the
line segment $1/2 - E_2$ lies in the interior of the triangle $\<0,\;
1/4,\; 1/5\>$.  Therefore there is no solution in this case.

When $E_2 = \<1/3,\; 2/5\>$, we have the following calculation.
$$
\align
& y_2  = \frac 13 + \frac 12 (x - \frac 23) = \frac x2 \\
& y_3 = \frac 1{l_3} x \qquad \text {and }\\
& y_3 = -y_1 -y_2 = \frac 12 - y_2 = \frac 12 - \frac x2 \qquad \text{hence}\\
& x = \frac {l_3}{l_3 + 2}  \\
& \bu = \frac 1 { 1-x} = \frac{l_3 + 2}{2}  
\endalign
$$
Since $v_2 \in \< 1/3,\; 2/5 \>$ and is not a vertex, we have $3 < \bu
< 5$, so the above gives $4 < l_3 < 8$..  Since $K(-1/2,\; 2/5,\;
1/s_3)$ is a knot, $s_3$ must be odd.  Therefore the only possibility
is that $s_3 = 7$, in which case $l_3 = 6$, $\bu = 4$, and the knot is
$K(-1/2,\; 2/5,\; 1/7)$.

When $E_2 = \<1/3,\; 2/7\>$, $1/2 - E_2$ lies on the edges $\< 0,\;
1/5 \>$ and $\< 1/5,\; 2/9 \>$.  Since $q_3 \leq 3$, we must have $E_3
= \< 0,\; 1/5 \>$.  Note that $l_2 = l_3$ and the slopes of these
edges add up to zero.  Therefore by Lemma 2.11 all solution surfaces
on these edges have the same boundary slope.  We have $K = K(-1/2,\;
1/5,\; 2/7)$, and $3<\bu<5$.

We have $s_2 \neq 6$ because $s_2$ is coprime with $q_2$.  Thus it
remains to consider the case that $E_2 = \<1/3,\; r_2/s_2\>$ for some
$s_2 \geq 8$.  It is clear from Figure 2.1 that $1/2 - E_2$ does not
intersect the interior of $\< 0,\; 1/s \>$ for $s \leq 4$.  Therefore we
may assume that $E_3 = \< 0,\; 1/s_3 \>$ for some $s_3 \geq 5$.  We now
have $l_2 = s_2 - q_2 \geq 5$ and $l_3 = s_3 - q_3 \geq 4$.  Since $\bu
> q_2 = 3$, by Lemma 5.1 we have
$$
\align
\bar e & \leq (1- \frac {q_2}{l_2} - \frac {q_3}{l_3}) + (- \frac 12 + \frac
1{l_2} + \frac 1{l_3}) \bu  \\
& = (1 - \frac 3{l_2} - \frac {1}{l_3}) + 3 (- \frac 12 + \frac
1{l_2} + \frac 1{l_3}) + (- \frac 12 + \frac 1{l_2} + \frac 1{l_3})
(\bu-3) \\
& = (-\frac 12 + \frac 2{l_3}) + (- \frac 12 + \frac 1{l_2} + \frac
1{l_3}) (\bu - 3)  \\ 
& \leq (-\frac 12 + \frac 24) + (-\frac 12 + \frac 15 + \frac 14)(\bu-3)
< 0
\endalign
$$
This contradicts Theorem 2.8.  Therefore there is no solution in this
case.
\endproof

\proclaim{Lemma 5.4} Suppose $v_1 \in L(-\frac 12)$, $\a_i>0$ for
$i=2,3$, and $q_3 \leq q_2 = 2$.  Then $K = K(-\frac 12,\; \frac
13,\; \frac 17)$, and $\bu = \frac 52$.  \endproclaim

\proof The minimum value of $y_2$ on edges of type $\< 1/2,\; r_2/s_2 \>$ is
achived at the vertex $\< 1/3 \>$ on $\<\frac 12, \frac 13\>$.  If
$q_3 = 2$ then $y_2 + y_3 \geq 2/3 > y_1$, so there is no solution to
$\sum y_i = 0$.  Therefore we must have $q_3 = 1$.  By the remark at
the beginning of the section we have $y_2 < 1/2$, so $E_2 = \< 1/2,\;
r_2/s_2 \>$ for some $r_2/s_2 < 1/2$, and $E_3 = \< 0,\; 1/s_3 \>$.

We have the following calculation.
$$
\align
& y_2  = \frac 12 - \frac 1{l_2} (x - \frac 12) \\
& y_3 = \frac 1{l_3} x \qquad \text {and}\\
& y_3 = \frac 12 - y_2 = \frac 1{l_2} (x - \frac 12), \qquad \text{hence}\\
& x = \frac {l_3}{2(l_3 - l_2)}  \\
& \bu = \frac 1 { 1-x} = \frac{2(l_3 - l_2)}{l_3 - 2l_2} = 2 + \frac
    {2l_2}{l_3 - 2l_2}  \\
& \a_2 = \frac {\bu-q_2}{l_2} = \frac{\bu-2}{l_2} = \frac 2{l_3 - 2l_2}  \\
& \a_3 = \frac {\bu-q_3}{l_3} = \frac{\bu-1}{l_3} = \frac 1{l_3 - 2l_2}  \\
& 0 \leq \bar e \leq e(\bu) = 1 - \frac 12 \bu + \a_2 + \a_3 =  \frac
{3 - l_2}{l_3 - 2l_2} 
\endalign
$$
Since $v_2 \in \< 1/2,\; r_2/s_2 \>$ and is not a vertex, we have $\bu < s_2 =
l_2 + 2$.  By the above formula of $\bu$, this gives 
$2+2l_2/(l_3 - 2l_2) < 2+ l_2$, hence
$$ l_3 - 2l_2 > 2  \tag *
$$
Note that the slope of $E_2$ is negative and the slope of $E_3$ is
positive, so by Lemma 2.11 the boundary slope $\delta$ of the surface
satisfies
$$
\delta \equiv 2(e_- - e_+) \equiv 2(\b_3 - \b_2) \equiv 2(\a_2 -
\a_3) = \frac 2{l_3 - 2l_2} \not\equiv 0  \qquad \text{mod $1$} 
$$
By Theorem 2.8 this means that $\frac 12 \leq \bar e < 1$, hence
$$ e(\bu) = \frac {3-l_2}{l_3 - 2l_2} = \frac 12 \times \frac
{3-l_2}{(l_3/2) - l_2} \geq \bar e \geq \frac 12$$
By (*) we have $(l_3/2)-l_2 > 0$, so the above inequality gives 
$$
l_3 \leq 6$$
Together with (*), this implies that $(l_2, l_3) =
(1,5)$ or $(1,6)$.  By definition $e(\bu) - \bar e$ equals the number
of full edges in $\cup e(\r_i)$.  In both cases $e(\bu) < 1$, so there
are no full edges.  In the first case we have $K = K(-1/2,\; 1/3,\;
1/6)$, which is not a knot.  In the second case we have $\bu=2.5$,
$s_2 = 3$, and $s_3 = 7$, hence the knot is $K = K(-1/2,\; 1/3,\;
1/7)$.  This solution gives the well-known $37/2$ toroidal surgery on
$K(-1/2,\; 1/3,\; 1/7)$.  \endproof

\proclaim{Lemma 5.5} Suppose $v_1 \in L(-\frac 12)$, $\a_i> 0$, 
and $E_i = \<0,\; 1/s_i\>$ for $i=2,3$.  Then 
$$
\align 
& K = K(-1/(2+1/n),\; 1/3,\; 1/3)  \qquad  \text{ $n$ odd, $n\neq-1$, $\bu =
2$; or }  \\
& K = K(-1/2,\; 1/3,\; 1/(3+1/n))  \qquad  \text{ $n$ even, $n \neq
0$, $\bu = 2$}
\endalign
$$
\endproclaim

\proof
Since $v_i$ is in the interior of the edge $E_i$ above, we have $y_i <
1/s_i$, so $\sum y_i = 0$ has no solution if $s_i \geq 4$ for
$i=2,3$.  Also since $v_1 \in L(-\frac 12)$, we have $\bu \geq 2$,
hence $s_i > 2$.  Therefore we may assume that $s_2 = 3$ and $s_3 \geq
3$.  We have
$$
\align
& y_2 + y_3 = \frac x{l_2} + \frac x{l_3} = (\frac 1{l_2} + \frac
1{l_3}) x = \frac 12 = - y_1 \\
& x = \frac {l_2 l_3}{2(l_2 + l_3)} \\
& \bu = \frac 1{1-x} = \frac {2(l_2 + l_3)} {2(l_2 + l_3) - l_2 l_3} \\
& \a_2 = \frac {\bu-q_2}{l_2} = \frac {l_3} {2(l_2 + l_3) - l_2 l_3} \\
& \a_3 = \frac {\bu-q_3}{l_3} = \frac {l_2} {2(l_2 + l_3) - l_2 l_3} \\
& \a_2 + \a_3 = \frac 12 \bu \\
& \bar e  \equiv e(\bu) = 1 - \frac 12 \bu + \a_2 + \a_3 = 1
\endalign
$$

Since both $E_2$ and $E_3$ have positive slope, $e_- \equiv \b_2+\b_3
\equiv -(\a_2+\a_3)$ mod $1$;
by Lemmas 2.11 the boundary slope of the surface satisfies $\delta
\equiv 2(e_- - e_+)
\equiv -2(\a_2 + \a_3) = -\bu$.  Since $\bar e \equiv 0$ mod $1$, by
Theorem 2.8 $\delta = -\bu$ must be an integer slope, and $\bar e = 0$.
Since $v_1 \in L(-\frac 12)$ and $v_2$ is in the
interior of $\< 0,\; 1/3 \>$, we have $2 \leq \bu < 3$, hence $\bu=2$.
From the formula of $\bu$ above we have $l_2 = l_3 = 2$.  Also
$$0 =\bar e = e(\bu) - \sum (|\r_i| - \b_i) = 1 - \sum(|\r_i| - \b_i)
$$
Hence there is an extra edge, which may end at either $\<\frac
12\>$ or $\<\frac 13\>$.  Therefore $K$ is one of the following knots.
$$
\align 
&K(-\frac 1{2+1/n},\; \frac 13,\; \frac 13)  \qquad  \text{
$n$ odd, and $n\neq -1$;}  \\ 
&K(-\frac 12,\; \frac 13,\; \frac 1{3+1/n})  \qquad  \text{ $n$
even, $n \neq 0$} 
\endalign
$$
The extra conditions on $n$ is to guarantee that $\r_i$ are allowable,
and $K$ is a knot.
\endproof

\proclaim{Proposition 5.6}  If $v_1 \in L(-\frac 12)$ and $\a_i > 0$ for
$i=2,3$, then $K$ is one of the knots in Lemma 5.3, 5.4 or 5.5.
\endproclaim

\proof
Because of symmetry we may assume $q_2 \geq q_3 \geq 1$.  By Lemma 5.2
we have $q_2 \leq 3$, so $q_2 = 3$, $2$ or $1$, which are covered by
Lemmas 5.3, 5.4 and 5.5, respectively.
\endproof

\head 6. No horizontal edges \endhead

In this section, we study the case that no $v_i$ is on a horizontal
edge.  As before, let $E_i = \<p_i/q_i,\; r_i/s_i\>$.  
By Lemma 2.10(4)
we may assume that $v_1$ is in the interior of a non-horizontal edge
in the graph $G$ shown in Figure 2.4, hence $1< \bu < 3$ and $\bu \neq
2$.  Since $q_i \leq \bu$, we have $q_i \leq 2$ for all $i$.

Similar to the previous sections, we define
$$
\align
\b_i(u) & = \frac {s_i - u}{s_i - q_i} \\
\a_i(u) & = \frac {u - q_i}{s_i - q_i} \\
e(u) & = 4-u - \sum \b_i(u) = 1 - u + \sum \a_i(u)
\endalign
$$
Since no $v_i$ is in the interior of a horizontal edge, by Definition
2.5 we have 
$$ 
\align
\bar e & = \sum e(\r_i) = (4-\bu) - \sum |\r_i| \\
& \leq 4-\bu - \sum \b_i = 1-\bu + \sum \a_i = e(\bu) \tag 6.1
\endalign
$$
Note also that $e(\bu) - \bar e$ is a nonnegative integer, which
equals the number of full edges in $\cup \r_i$.  

Let $\delta$ be the boundary slope of $F(\r_1, \r_2, \r_3)$.  Then 
$$
\delta \equiv 2(e_- - e_+) \equiv - 2 \sum \text{sign}(r_i/s_i -
p_i/q_i)\a_i \qquad \text{mod $1$} $$

\proclaim{Lemma 6.1}  Suppose $v_1 \in G - L$.  If $q_i = 2$, then 
$s_i = 3$.
\endproclaim

\proof
Assume to the contrary that $q_2 = 2$ and $s_2 > 3$.  We have $E_2 =
\<\pm 1/2,\; r_2/s_2\>$.  From Figure 2.1 we see that any point $(x,
y_2)$ in the interior of $E_2$ satisfies $x/2 < |y_2| < x$.

First assume that $s_3 = 3$.  Denote by $\pm Q = \<0,\; \pm \frac 13\>
\cup \< \pm \frac 12, \pm \frac 13\> \cup \< \pm \frac 12, \pm \frac
23\>$.  Note that any point $(x,y)$ on $\pm Q$ satisfies $x/2 \leq
|y| \leq x$.  

Since $\bu > q_2 = 2$, we have $s_1 = 3$, so $v_1 \in G-L$ implies that $v_1
\in -Q$.  We assumed $s_3=3$ and by Lemma 2.10 we have $|y_3| \leq \frac
23$, hence $v_3 \in \pm Q$.  Thus $x/2 \leq |y_i| \leq x$ for $i=1,3$.
This implies that if $y_3<0$ then $-y_1 -y_3 \geq (x/2) + (x/2) = x$,
and if $y_3 < 0$ then $|-y_1 - y_3| \leq x - x/2 = x/2$, either case
contradicting the fact that $y_2 = -y_1 - y_3$ satisfies $x/2 < |y_2|
< x$.

We may now assume that $s_3 > 3$.  Consider the function $e(u) = 1-u +
\sum \a_i(u)$, where $\a_i(u) = (u-q_i)/(s_i-q_i)$.  By (6.1)
we have $e(\bu) \geq \bar e \geq 0$.

When $u=2$ we have $\a_1(2) \leq \frac 12$, $\a_2(2) = 0$ because $q_2
= 2$, and $\a_3(2) \leq \frac 13$ because $s_3\geq 4$, hence $e(2) =
(1-2) + \sum \a_i(u) < 0$.  

Now calculate $e(3)$.  Since $s_2$ is coprime with $q_2=2$ and
$s_2>3$, we have $s_2 \geq 5$.  Hence $\a_2(3) \leq \frac 13$.  Also,
$\a_3(3) \leq \frac 23$ because $s_3>3$, and $\a_1(3) \leq 1$.  Hence
$e(3) = (1-3) + \sum \a_i(3) \leq 0$.  By the linearity of $e(u)$ we
have $\bar e \leq e(\bu) < 0$ because $2= q_2 <\bu< s_1 = 3$.  This 
contradicts Theorem 2.8.
\endproof

\proclaim{Lemma 6.2}  Suppose $v_1 \in \Int \< - \frac 12, -\frac
23 \>$, and $|y_1| > y_2 \geq y_3>0$.  Then $\bu = 2.5$, and $K =
K(-2/3,\; 1/3,\; 1/4)$.
\endproclaim

\proof
The equation for $E_1$ is $y_1 = -x$.  By Lemma 2.10(2) we have $|y_1|
+ |y_2| \leq 1$, hence $y_2
\leq 1 - x$.  Since $y_2 + y_3 = - y_1$ and $y_2 \geq y_3$, we have
$y_2 \geq \frac 12 (-y_1) = \frac 12 x$.  From Figure 2.1 we see that
$E_2 = \< \frac 12, \frac 13 \>$ or $\<
0, \frac 13 \>$.

If $E_2 = \< 0,\; \frac 13 \>$ then $y_3 = - y_1 - y_2 = x - x/2 = x/2 =
y_2$, so $E_3 = E_2$.  In this case $E_1 = - E_2 - E_3$, so there are
infinitely many solutions, all giving the same slope.  We have $\a_1 =
\bu - 2$, $\a_2 = \a_3 = (\bu-1)/2$, hence 
$$\bar e \leq e(\bu) = 1-\bu + \sum \a_i
= \bu -2 < 1$$  
Therefore there are no extra edges, hence $K=K(-\frac
23, \frac 13, \frac 13)$.  Since this is a link of two components.
there is no solution in this case.

Now assume $E_2 = \< \frac 12, \frac 13 \>$.  Then $y_3 = -y_1 - y_2 =
x - (1-x) = 2x -1 < \frac 13$ because $x < \frac 23$.  Hence from
Figure 2.1 we see that $E_3 = \< 0,\; 1/s_3 \>$ for some
$s_3 \geq 3$.  Define $l_i = s_i - q_i$.  We have the following
calculation. 
$$
\align
& -y_1 - y_2 = 2x-1 = \frac{x}{l_3}  = y_3 \\
& x = \frac 1{2 - 1/l_3} = \frac {l_3}{2l_3 -1} \\
& \bu = \frac 1{1-x} = 2 + \frac 1{l_3 - 1} \\
& \a_1 = \frac {\bu-2}{l_1} = \frac 1{l_3 - 1} \\
& \a_2 = \frac {\bu-2}{l_2} = \frac 1{l_3 - 1} \\
& \a_3 = \frac {\bu-1}{l_3} = \frac 1{l_3 - 1}  \\
& 0 \leq \bar e \leq e(\bu) = 1-\bu + \sum \a_i = -1 + \frac 2{l_3 - 1}
\endalign
$$
This gives $l_3 = 2$ or $3$.  When $l_3 = 2$, $\bu = 3$, which is a
contradiction because $\bu <3$.  When $l_3 = 3$, we have $\bu = 2.5$
and $\bar e = 0$, so there is no extra edge, hence 
$K = K(r_1/s_2,\; r_2/s_2,\;, r_3/s_3) = K(-2/3,\; 1/3,\; 1/4)$.
\endproof

\proclaim{Lemma 6.3}  Suppose $v_1 \in \Int \<-1/2,\; -1/3 \>$, and
$|y_1| > y_2 \geq y_3>0$.  Then $K = K(-1/3,\; 1/3,\; 1/7)$ and $\bu =
2.5$.
\endproclaim

\proof
Since $-y_1 > y_2 \geq y_3 > 0$, $E_2$ and $E_3$ must be below the
edge $E_1 = \<\frac 12, \frac 13\>$, hence we must have
$E_i = \< 0,\; 1/s_i \>$ for $i=2,3$.  We have the following
calculation.
$$
\align
& \sum y_i = (x-1) + x/l_2 + x/l_3 = 0  \\
& x = \frac 1{(1/l_2) + (1/l_3) + 1}    \\
& \bu = \frac 1{1-x} = 1 + \frac{l_2 l_3}{l_2 + l_3} \\
& \a_1 = \frac {\bu-2}{l_1} = \bu - 2 =  \frac{l_2 l_3}{l_2 + l_3} - 1\\
& \a_2 = \frac {\bu-1}{l_2} = \frac {l_3}{l_2 + l_3} \\
& \a_3 = \frac {\bu-1}{l_3} = \frac {l_2}{l_2 + l_3} \\
& \delta \equiv 2(- \a_1 - \a_2 - \a_3) \equiv \frac {- 2 l_2 l_3}{l_2 +
l_3} \equiv -2\bu \\ 
& 0 \leq \bar e \leq  1 - \bu + \sum \a_i = 0
\endalign
$$
The last inequality gives $\bar e = 0$, so by Theorem 2.8 the boundary
slope of the surface must be an integer, hence by Lemma 2.11 we have
$\delta \equiv 2\bu \equiv 0$ mod $1$.  Since $2<\bu<3$, we have $\bu
= 2.5$.  The only solutions for $ \bu = 1 + l_2 l_3/(l_2 + l_3) = 2.5
$ are $(l_2, l_3) = (2,6)$ or $(3,3)$.  Since $e=0$, there is no extra
edge, so in the second case we would have $K = K(-1/3,\; 1/4,\; 1/4)$,
which is not a knot.  Therefore the only solution in this case is $K =
K(-1/3,\; 1/3,\; 1/7)$ and $\bu = 5/2$.
\endproof

\proclaim{Lemma 6.4}  Suppose $v_1 \in G$, and $|y_1| > y_2 \geq
y_3 >0$.  Then $v_1 \notin \Int \< 0,\; -1/2 \>$.  
\endproclaim

\proof
If $v_1 \in \Int \< 0,\; -1/2 \>$, then $|y_1| > y_i > 0$ implies that
$E_i = \<0, 1/s_i \>$ for some $s_i \geq 3$.  We have
$y_1 = -x$, and $y_i = x/l_i$ for $i=2,3$, hence from the equations $y_1 +
y_2 + y_3 = 0$ and $x>0$ we have
$$
-1 + \frac 1{l_2} + \frac 1{l_3} = 0
$$
Since $s_i \geq 3$, this gives $l_2 = l_3 = 2$.  Note that $\a_1 =
\bu-1$, and $\a_2 = \a_3 = (\bu - 1)/2 < 1/2$, hence 
$$
\bar e = 4-\bu -\sum |\r_i| \equiv 1 - \bu + \sum \a_i = \a_2 + \a_3 =
\bu -1 \not\equiv 0
\qquad \text{mod $1$}
$$
By Theorem 2.8 this implies that the boundary slope $\delta$ of the
surface is not an integer slope.  On the other hand, since $E_1$ has
positive slope and $E_2, E_3$ have negative slope, by Lemma 2.11 we
have 
$$\delta \equiv 2(e_- - e_+)
= 2(\b_2 + \b_3 - \b_1) \equiv 2(-\a_2 - \a_3 + \a_1) = 0 \qquad
\text{mod $1$}
$$
so $\delta$ is an integer slope, which is a contradiction.
\endproof

\proclaim{Lemma 6.5}  Suppose $v_1 \in \Int \< 0,\; -1/3 \>$, and
$|y_1| > y_2 \geq y_3 > 0$.  Then $K = K(-1/3,\; 1/4,\; 1/7)$ or
$K(-1/3,\; 1/5,\; 1/5)$, and the boundary slopes are the same as the
pretzel slopes corresponding to the candidate systems in Proposition
2.9. 
\endproclaim

\proof
Similar to the proof of Lemma 6.4, we have
$$
\sum y_i = -\frac 12 + \frac 1{l_2} + \frac 1{l_3} = 0
$$ 
which gives $(l_2, l_3) = (3,6)$ or $(4,4)$.  We have
$$
0 \leq \bar e \leq e(\bu) = 1 - \bu + \sum \a_i(\bu) = 1 - \bu
+ \frac {\bu-1}{2} + \frac {\bu-1}{l_2} + \frac {\bu-1}{l_3} = 0
$$
Therefore there are no extra edges.  The knots are $K = K(-1/3,\;
1/4,\; 1/7)$ and $K(-1/3,\; 1/5,\; 1/5)$.  The boundary slopes are the
same for all $\bu$, which is also the boundary slope of the candidate
system at $\bu =1$, as given in Proposition 2.9.
\endproof

\proclaim{Lemma 6.6}  Suppose $v_1 \in G$, and $|y_1| > y_2 \geq
y_3 >0$.  Then $v_1 \notin \Int \< -1, -\frac 12 \>$.
\endproclaim

\proof
If $v_1 \in \Int \< -1, -\frac 12 \>$ then 
$1<\bu < 2$, and we have $E_i = \< 0,\; 1/s_i \>$ for
$i=2,3$.  If both $s_i > 2$ then it is easy to show that $|y_1| >
|y_2| + |y_3|$ for $0 < x < \frac 12$, so there is no solution to
$\sum y_i = 0$ when $1<\bu < 2$.  Hence we must have $s_2 = 2$, so
$y_1 = -1+x$, $y_2 = x$, and $y_3 = x/l_3$.  We have the following
calculations.  
$$
\align
& \sum y_i = (-1+x) + x + \frac x{l_3} = 0 \\
& x = \frac {l_3}{2l_3 + 1} \\
& \bu = \frac 1{1-x} = 2 - \frac 1{l_3 + 1} = 2 - \frac 1{s_3} \\
& \a_1 = \frac {\bu-q_1}{l_1} = \frac {(2-1/s_3)-1}{1} = 1 - 
\frac {1}{s_3}  \\
& \a_2 = \frac {\bu-q_2}{l_2} = \frac {(2-1/s_3)-1}{1} = 1 - 
\frac {1}{s_3}  \\
& \a_3 = \frac {\bu-q_3}{l_3} = \frac {(2-1/s_3)-1}{s_3 - 1} = 
\frac {1}{s_3}   \\
& \delta \equiv 2(- \a_1 - \a_2 - \a_3) =  \frac {2}{s_3} \\
& \bar e \equiv e(\bu) = 1 - \bu + \sum \a_i = 1 \equiv 0 \qquad
\text{mod $1$}  
\endalign
$$
By Theorem 2.8, $\bar e \equiv 0$ mod $1$ implies that $\delta$ is an
integer slope, hence from $\delta \equiv 2/s_3$ mod $1$ and $s_3 \geq
2$ we see that $s_3 = 2$.  Since $\bar e = 1$ there is one extra edge,
but since $(r_1/s_1,\; r_2/s_2,\; r_3/s_3) = (-1/2,\; 1/2,\; 1/2)$ or
$(-1/2,\; 1/2,\; 1/4)$, adding one extra edge will make a link of type
$K(-1/2,\; 1/2,\; 1/(2+1/n))$ or $K(-1/2,\; 1/2,\; 1/(4+1/n))$, which
has at least two components.  Therefore there is no solution in this
case.
\endproof

\proclaim{Proposition 6.7}  Suppose some $v_i \in G-L$, and $\bu > 1$.
Then $K$ and $\bu$ are equivalent to one of the pairs in Lemmas 6.2,
6.3 and 6.5.
\endproclaim

\proof
Since $\bu > 1$, we must have $y_j \neq 0$.  Up to permutation and
changing of signs of the parameters of $K$ we may assume that $- y_1
\geq |y_i|$ for $i=2,3$.  Before this modification we have some $v_i
\in G-L$.  We need to show that $v_1 \in G-L$ after the modification.

The assumption of $v_i \in G-L$ implies that $\bu < 3$, hence $p_1/q_1
\in \{0, -\frac 12, -1\}$.  By Lemma 2.10(2) we have $|y_1|
\leq \frac 23$, so from Figure 2.1 we see that if $p_1/q_1 = -1$
then $E_1 = \< -1, -\frac 12 \>$, hence $v_1 \in G-L$.  If
$p_1/q_1 = -\frac 12$ then by Lemma 6.1 we have $r_1/s_1 = -\frac 13$
or $-\frac 23$, so again $v_1 \in G-L$.  If $p_1/q_1 = 0$ and
$v_1 \notin G-L$ then from Figure 2.1 and the assumption of $|y_i| \leq
|y_1|$ we see that each $v_j$ is in $\< 0, \pm \frac 1{s_j} \>$
for some $s_j \geq 4$, which contradicts the assumption that some
$v_i$ is in $G-L$ before the modification of parameters of $K$.

We now have $v_1 \in G-L$, and $-y_1 \geq |y_i|$ for $i=2,3$.  Since
$\bu > 1$, no $v_i$ is on the horizontal line $L(0)$ as otherwise
the corresponding parameter of $K$ would be $0$, contradicting the
assumption that the length of $K$ is 3.  Hence from $\sum y_i =0$ we
have that $|y_1| > y_i > 0$.  Permuting the second and third
parameters of $K$ if necessary we may assume $y_2 \geq y_3$.
Therefore we have $|y_1| > y_2 \geq y_3 > 0$.  

There are only 5 edges in $G-L$, so $E_1$ must be one of them, which
have been discussed in Lemmas 6.2 -- 6.6 respectively.  The above
discussion shows that $y_i$ satisfy the conditions of the lemmas,
hence the result follows from these lemmas.
\endproof

\head 7.  The classification \endhead

\proclaim{Lemma 7.1} Let $K$ be a Montesinos knot of length $3$.  Then
$E(K)$ contains a candidate surface $F$ of genus one with boundary
slope $\delta$ if and only if $(K, \delta)$ is equivalent to one of
the pairs listed in Theorem 1.1.
\endproclaim

\proof
When $\bu \leq 1$, the knots are given in Proposition 2.9.  When $\bu
> 1$, by Lemma 2.10(4) we may assume that $v_1 \in G$.  Propositions
3.7, 4.4 and 5.6 covered the case of $v_1 \in L$, and Proposition 6.7
covered the case of $v_1 \in G-L$.  The list in Theorem 1.1 contains
all the knots given in these propositions.  Here are more details.

Parts (1) and (2) of Theorem 1.1 include the knots in Proposition 2.9,
as well as $K(-1/3,\; 1/5,\; 1/5)$ and $K(-1/3,\; 1/4,\; 1/7)$ in
Proposition 4.4 and Lemma 6.5.  Note that the boundary slopes for the
last two are the same as those in the list, but the $\bu$ values are
different, which is allowed by the remark before the statement of
Theorem 1.1.

Part (3) is given by Lemma 3.2, and part (4) by Lemma 3.5.  (5), (6)
and (9) are in Lemma 3.6.  (7) and (8) are from Lemma 5.5, (10) from
Lemma 5.3, (11) from Lemma 5.4, (12) from Lemma 6.2, and (13) from
Lemma 6.3.  Note that the knot $K(-1/2,\; 1/5,\; 2/7)$ in Lemmas 3.6
and 5.3 is included in (5) (with $n=2$) because they all have the same
boundary slope.  The boundary slopes can be calculated using the
algorithm of Hatcher and Oertel in Lemma 2.11.
\endproof

\demo{Proof of Theorem 1.1} By Lemma 7.1 we need only show that the
candidate surfaces in Lemma 7.1 are incompressible in $E(K)$ when $K$
is hyperbolic.  Since the candidate systems are already determined in
the proofs of the lemmas and the propositions above, it is straight
forward to follow the procedure of [HO] to verify the
incompressibility of the candidate surface $F$.  For each individual
knot the toroidal slopes can also be verified using the computer
program of Nathan Dunfield [Dn].

Here are more details.  If $\bu \geq q_i$ for some $i$ then the path
$\gamma_i$ is a constant path, so it follows from [HO, Proposition
2.1] that $F$ is incompressible.  This takes care of all the cases
except (1), (2), (7), (12) and (13) in Theorem 1.1.  Recall from [HO,
Page 463] that the $r$-value of a path is determined by this rule:
Extende the last segment of the path to meet the right-hand border of
the the diagram in a point whose slope has denominator $r$.  (For
length 3 knot we do not have vertical line segments in $\gamma_i$ and
hence $r\neq 0$.)  By [HO, Corollary 2.4] $F$ is incompressible unless
the $r$-value cycle is $(1,1,r_3)$ or $(1,2,r_3)$.  In cases (1) and
(2) the $r$-value cycle is $(|q_1|-1,\; |q_2|-1,\; |q_3|-1)$, so $F$
is incompressible unless $|q_1| = 2$ and $|q_2| = 3$.  ($|q_2|\neq 2$
since $K$ is a knot.)  By [BO, Proposition 2.7] $F$ is incompressible
unless the slopes of $\gamma_1$ is of opposite sign to those of
$\gamma_2$ and $\gamma_3$ and $r_3 = |q_3|-1 = 2$ or $4$.  It follows
that $F$ is incompressible unless $K$ is the non-hyperbolic knot
$K(-1/2,\, 1/3,\, 1/3)$ or $K(-1/2,\, 1/3,\, 1/5)$, which has been
excluded since we assumed that $K$ is hyperbolic.  In case (7) $r_i
\neq 1$ unless $K = K(-1/3,\, 1/3,\, 1/3)$, in which case the
$r$-cycle is $(1,2,2)$.  Since the slopes of the three edges are all
positive (the first one is on $\< -1/2,\; -1/3 \>$ and the other two
are on $(\<0/1,\;, 1/3\>$), by [HO, Proposition 2.7] $F$ is
incompressible.  In case (12) the $r$-cycle is $(1,1,3)$ and the first
two paths are edges of the same slope $-1$, so $F$ is incompressible
by [HO, Proposition 2.6].  In case (13) the $r$-cycle is $(1,2,6)$ and
the edges all have positive slopes, hence again the incompressibility
of $F$ follows from [HO, Proposition 2.7].
\endproof

If $F$ is a surface in a 3-manifold $M$, denote by $M|F$ the manifold
obtained by cutting $M$ along $F$.  Similarly, if $C$ is a set of
curves on a surface $F$ then $F|C$ denotes the surface obtained by
cutting $F$ along $C$.  All surfaces in 3-manifolds below are assumed
compact, connected, orientable, and properly embedded.  Recall that a
surface in $M$ is {\it essential\/} if it is incompressible,
$\bdd$-incompressible, and is not boundary parallel.  Denote by $|\bdd
F|$ the number of boundary components of $F$.

\proclaim{Lemma 7.2}  Let $M$ be a compact orientable 3-manifold with
$\bdd M = T$ a torus.  Let $F$ be a genus one separating essential
surface in $M$ such that $|\bdd F| \leq 4$.  Let $M_1, M_2$ be the
components of $M|F$, and let $A$ be a component of $\bdd M_1 - F$.  If
$F\cup A$ is incompressible in $M_1$, then $M$ contains a closed
essential surface.
\endproclaim

\proof This is due to Gordon, and is true for any number of components
on $\bdd F$.  If $|\bdd F| = 2$ then by assumption $F\cup A$ is
incompressible and the result follows, so we assume $|\bdd F| = 4$.
Let $A'$ be the annulus on $T$ which contains $A$ and such that $\bdd
A' = \bdd F - \bdd A$.  Then we can push the part of $F \cup A$ near
$A$ into the interior of $M$ to obtain a surface $F'$ with $\bdd F' =
A'$, and then push $F' \cup A'$ into the interior of $M$ to obtain a
closed surface $F''$.  Number the annuli $T|\bdd F$ successively as
$A_1, ..., A_4$, with $A = A_1$.  One can show that $M|F''$ consists
of two components $M''_1, M''_2$, such that $M''_1$ is obtained by
gluing $T\times I$ to $M_1$ along the annulus $A_3$ and a nontrivial
annulus on $T\times 0$, and $M''_2$ is obtained by gluing a $A' \times
I$ to $M_2$, where $A'$ is an annulus, $A'\times 0$ identified to
$A_2$, and $A'\times 1$ to $A_4$.  An innermost circle argument shows
that if $F\cup A$ is incompressible then $F''$ is incompressible in
both $M''_i$, hence is incompressible in $M$.  \endproof

\proclaim{Lemma 7.3}  Let $M$ be a compact orientable irreducible
3-manifold with $\bdd M = T$ a torus.  Let $F$ be a genus one
separating incompressible surface in $M$ with boundary slope $\delta$,
and let $\hat F$ be the corresponding torus in the Dehn filling
manifold $M(\delta)$.  If (i) $M$ contains no closed incompressible
surface, and (ii) $F$ has at most four boundary components, then $\hat
F$ is incompressible in $M(\delta)$.
\endproclaim

\proof  
Since $F$ is separating, $|\bdd F| =2$ or $4$.  We assume the latter,
as the proof for the former case is similar and simpler.
Let $M_1, M_2$ be the components of $M|F$, and let $A_1, ..., A_4$ be
the annuli $T|\bdd F$, labeled so that $\bdd M_1 = F \cup A_1 \cup
A_3$.  

Since $M$ contains no closed essential surface, each $M_i$ is a
handlebody of genus 3, and by Lemma 7.2 the surface $F_1 = F \cup A_1$
is compressible in $M_1$.  Let $D$ be a compressing disk of $F_1$.
If $D$ is separating then since $M_1$ is a handlebody and $\bdd M_1 -
F_1 = A_3$ is connected, we can find a nonseparating compressing disk
in a component of $M_1|D$ disjoint from $A_3$.  Therefore we may
assume without loss of generality that $D$ is a nonseparating
compressing disk.  It follows that after attaching a 2-handle to $M_1$
along $A_3$ the resulting manifold $M'$ has compressible boundary,
because $D$ remains a compressing disk of $\bdd M'$ in $M'$.

We have $A_1 \subset \bdd M'$.  We want to show that $F_2 = \bdd M' -
A_1$ is incompressible in $M'$.  Consider the surface $F_3 = F\cup A_3 =
\bdd M_1 - A_1$.  For the same reason as above, we know
that $F_3$ is compressible in $M_1$.  By assumption $F_3 - A_3 = F$ is
incompressible.  Therefore by the Handle Addition Lemma (see [Ja] or
[CG]), we know that after attaching a 2-handle to $M_1$ along $A_3$,
the resulting surface $F_2 = \bdd M' - A_1$ is incompressible in $M'$.

We have shown that $\bdd M'$ is compressible, and $\bdd M' - A_1 =
F_2$ is incompressible.  Applying the Handle Addition Lemma again, we
see that after attaching a 2-handle to $A_1$ the boundary of the
resulting manifold $M''$ is incompressible.  Note that $M''$ is a
component of $M(\delta)|\hat F$, and $\bdd M'' = \hat F$.  For the
same reason, $\hat F$ is incompressible in the other component of
$M(\delta)| \hat F$.  It follows that $\hat F$ is incompressible in
$M(\delta)$.
\endproof

{\it Remark.\/} The above result is similar to a special case of
Proposition 2.2.1 of [CGLS].  However, that proposition requires that
the number of boundary components of the surface is minimal among all
incompressible surfaces with the same boundary slope.  In our case
there is no guarantee that there is no higher genus surface with fewer
boundary components of the same slope.  Lemma 7.3 is probably false if
there is no constraint about the number of components in $\bdd F$.

\demo{Proof of Theorem 1.2}
If $K_{\delta}$ is toroidal then clearly there is a toroidal
incompressible surface in the exterior of $K$, so by Theorem 1.1 the
pair $(K, \delta)$ must be one of those in the list.  

To prove the other direction, we would like to show that if $(K,
\delta)$ is in the list of Theorem 1.1 then the corresponding toroidal
incompressible surface $F = F(\r_1, \r_2, \r_3)$ gives rise to an
incompressible torus $\hat F$ in $K_{\delta}$.  By Oertel [Oe], the
exterior of $K(t_1, t_2, t_3)$ contains no closed essential surface.
Therefore by Lemma 7.3 it suffices to show that $F$ has at most four
boundary components.

Let $m_i$ be the number defined before the statement of Lemma 2.4, and
let $n=\lcm(m_1, m_2, m_3)$.  Then by Lemma 2.4(3) we have $|\bdd F|
\leq 2n$.  Therefore if $n\leq 2$ then by Lemma 7.3 the surface $\hat
F$ is an incompressible torus in $K_{\delta}$ and we are done.  By
definition $m_i$ can be easily calculated from $\bu$ and $E_i =
\<p_i/q_i,\; r_i/s_i\>$, which can be found in the proof of the
corresponding lemma for that $(K, \delta)$.  We leave it to the
reader to check that $m_i \leq 2$ for all $i$ in all the cases listed
in the theorem, except that $m_3 = 4$ in case (13).  (For each
individual knot, one may also use Dunfield's program [Dn] to calculate
$n = \lcm(m_1, m_2, m_3)$, which is shown as ``number of sheets'' in
the program.)  Therefore Theorem 1.2 follows from Lemma 7.3, except in
Case (13) of Theorem 1.1.

For Case (13), let $F' = F'(\r_1, \r_2, \r_3)$ be the surface in the
exterior of $K$ constructed using the candidate system $(\r_1, \r_2,
\r_3)$ given by Lemma 6.3, such that $|\bdd F'| = n = 4$.  By the
proof of Lemma 2.4 we have $F = F'$, or its double cover if $F'$ is
nonorientable.  In the first case we have $|\bdd F| =4$ and the result
follows from Lemma 7.3.  Hence we assume that $\hat F$ is a double
cover of $\hat F'$, so $\hat F'$ is a Klein bottle in $M =
K_{\delta}$.  On the other hand, from Theorem 1.1 we see that in this
case $\delta = 1$, hence $H_1(K_{\delta}, \Bbb Z_2) = 0$, which is a
contradiction because by duality a $\Bbb Z_2$-homology sphere cannot
contain a Klein bottle.  \endproof

\enddocument